\documentclass[
]{amsart}
\usepackage{algorithm, algorithmic}
\usepackage[latin1]{inputenc}
\usepackage{amsmath}\usepackage{amsthm}\usepackage{amssymb}\usepackage{mathrsfs}\usepackage{amscd}\usepackage{graphicx}\usepackage{subfigure}\usepackage{amsfonts}\usepackage{amsxtra}\usepackage{color}
\usepackage{tikz}
\usetikzlibrary{matrix}
\usepackage[american]{babel}

\usepackage[left=3.1cm, right=3.1cm,top=3.1cm,bottom=3.1cm]{geometry}

\DeclareMathOperator{\rank}{rank}

\DeclareMathOperator{\Hom}{Hom}

\DeclareMathOperator{\sym}{sym}

\DeclareMathOperator{\FFP}{FFP}

\DeclareMathOperator{\Null}{null}

\DeclareMathOperator{\diag}{diag}\DeclareMathOperator*{\spann}{span}

\DeclareMathOperator{\Pol}{Pol}

\DeclareMathOperator{\GL}{GL}
\DeclareMathOperator{\Tr}{Tr}

\newcommand{\R}{\mathbb{R}}
\newcommand{\C}{\mathbb{C}}
\newcommand{\N}{\mathbb{N}}

\newcommand{\G}{\mathcal{G}}

\renewcommand{\d}{\mathrm{d}}

\theoremstyle{definition}
\newtheorem{definition}{Definition}[section]
\newtheorem{remark}[definition]{Remark}
\theoremstyle{plain}\newtheorem{theorem}[definition]{Theorem}\newtheorem{lemma}[definition]{Lemma}\newtheorem{corollary}[definition]{Corollary}\newtheorem{proposition}[definition]{Proposition}

\begin{document}

\title[]{Reproducing kernels for the irreducible components of polynomial spaces on unions of Grassmannians}

%

\author[M.~Ehler]{Martin Ehler}
\address[M.~Ehler]{University of Vienna,
Department of Mathematics,
Oskar-Morgenstern-Platz 1, 
A-1090 Vienna, \textit{Phone: +43 1 4277 50729}, \textit{Fax: +43 1 4277 850729}
}
\email[Corresponding author]{martin.ehler@univie.ac.at}
\author[M.~Gr\"af]{Manuel Gr\"af}
\address[M.~Gr\"af]{Austrian Academy of Sciences,
Acoustics Research Institute, Wohllebengasse 12-14, 
A-1040 Vienna, Austria}
\email{mgraef@kfs.oeaw.ac.at}
\begin{abstract}


The decomposition of polynomial spaces on unions of Grassmannians $\mathcal G_{{k_1},d}\cup\ldots\cup \mathcal G_{{k_r},d}$ into irreducible orthogonally invariant subspaces and their reproducing kernels are investigated. We also generalize the concepts of cubature points and $t$-designs from single Grassmannians to unions. We derive their characterization as minimizers of a suitable energy potential to enable $t$-design constructions by numerical optimization. We also present new analytic families of $t$-designs for $t=1,2,3$. 

\end{abstract}

\maketitle

\section{Introduction}
Polynomial approximation from samples on manifolds and homogeneous spaces has already been extensively studied, cf.~\cite{deBoor:1993aa,DeVore:1993ab,Geller:2011fk,Maggioni:2008fk} and references therein. Constituting distinct sampling rules, the concepts of cubatures and designs have also been widely investigated, cf.~\cite{Delsarte:1977aa,Engels:1980uq,Filbir:2010aa,Hoggar:1982fk,Konig:1999fk,Neumaier:1988kl,Pesenson:2012fp}, where polynomials are integrated exactly by a finite sum over the sampling values. 
However, many open questions remain when dealing with polynomials on unions of non-connected manifolds. 

Orthogonal projectors with fixed rank are used in many applications for analysis and dimension reduction purposes, cf.~\cite{Harandi:2013wo,Veeraraghavan:2008aa}, leading to a function approximation problem on a single Grassmannian manifold. Projections with varying target dimensions are more flexible and may offer a larger range of applications. Therefore, we shall study unions of Grassmannians. 

By studying the structure of polynomial spaces on the union of Grassmannians, some of our findings generalize results in \cite{James:1974aa}. In particular, we shall verify that the multiplicities of the irreducible representations of the orthogonal group occurring in an orthogonally invariant reproducing kernel Hilbert space on unions of Grassmannians coincide with the ranks of the kernel's Fourier coefficients. This enables us to actually determine the multiplicities in the space of polynomials of degree $t$. Moreover, we construct the underlying reproducing kernels for the irreducible components.  While cubatures and $t$-designs in single Grassmannians have been studied in \cite{Bachoc:2005aa,Bachoc:2006aa,Bachoc:2004fk,Bachoc:2002aa}, we shall also investigate these concepts in unions of Grassmannians. We derive a characterization as minimizers of an energy functional induced by a reproducing kernel. By numerically minimizing the energy functional, we compute candidates for $t$-designs, i.e., $t$-designs up to machine precision. We are then able to
  describe these candidates analytically and check that the energy functional is
  minimized exactly.

It should be mentioned that the topic shares some common theme with Euclidean designs, cf.~\cite{Neumaier:1988kl}, where unions of spheres with varying radii in Euclidean space are considered, see also \cite{Bajnok:2006jt,Bajnok:2007wq,Bannai:2006la,Bannai:2012wd}. The structure of the polynomial spaces on unions of spheres have been investigated in \cite{Delsarte:1989jw}, but the ideas in those proofs do not work for unions of Grassmannians, whose structure appears to be more involved. 

The outline is as follows. 
In Section \ref{sec:2} we recall some facts on polynomial spaces on single
Grassmannians and their irreducible decompositions. Section \ref{sec:uni} is dedicated to some elementary results on polynomial spaces on unions of Grassmannians. Direct consequences of irreducible decompositions of
polynomials on symmetric matrices are studied in Section \ref{sec:H}. In Section \ref{sec:L2kernels} we determine the multiplicities of the polynomial spaces on unions and construct the underlying reproducing kernels for the irreducible components. In Section \ref{sec:cub} we introduce cubatures and $t$-designs on unions of Grassmannians and derive a characterization as minimizers of an energy functional induced by a reproducing kernel. We compute some analytical minimizers in Section \ref{sec:examples num}. 
\section{Polynomials on single Grassmannians}\label{sec:2}
This section is dedicated to summarize some facts about single Grassmannians, see, for instance, \cite{Bachoc:2002aa,James:1974aa}.   
The Grassmannian space of all $k$-dimensional linear subspaces of $\R^d$ is naturally
identified with the set of orthogonal projectors on $\R^d$ of rank $k$ denoted by
\[
\G_{k,d} := \{ P \in \R^{d\times d}_{\sym} \;:\; P^{2}=P ;\; \Tr(P)=k \},
\]
where $\R^{d\times d}_{\sym}$ is the set of symmetric matrices in $\R^{d \times d}$. Each Grassmannian $\mathcal{G}_{k,d}$ admits a unique 
orthogonally invariant probability measure $\sigma_{k,d}$ induced by the Haar
(probability) measure $\sigma_{\mathcal O(d)}$ on the orthogonal group
$\mathcal{O}(d)$, i.e., for any $Q \in \G_{k,d}$ and measurable function
$f$, we observe
\[
\int_{\G_{k,d}} f(P) \d\sigma_{k,d}(P) = \int_{\mathcal O(d)} f(O Q O^{\top}) \d
\sigma_{\mathcal O(d)}(O).
\]
The space of complex-valued, square-integrable functions $L^{2}(\G_{k,d})$,
endowed with the inner product $(f, g)_{\G_{k,d}}$, decomposes into orthogonally
invariant subspaces
\begin{equation}\label{eq:decomp L2 single}
L^{2}(\G_{k,d}) = \bigoplus_{\ell(\lambda) \le \min\{k,d-k\}}\!\!\!\!\!\!\! H_{\lambda}(\G_{k,d}),\qquad
H_{\lambda}(\G_{k,d}) \perp H_{\lambda'}(\G_{k,d}), \quad \lambda \ne \lambda',
\end{equation}
where $H_{\lambda}(\G_{k,d})$ is equivalent to $\mathcal H_{2\lambda}^{d}$, the
irreducible representation of $\mathcal O(d)$ associated to the partition $2\lambda=(2\lambda_{1},\dots,2\lambda_{t})$,
cf.~\cite{Bachoc:2002aa,James:1974aa}. Note that two representations are equivalent if there is a linear isomorphism that commutes with the group action. 
 A partition of $t$ is an integer vector $\lambda=(\lambda_1,\dots,\lambda_t)$
with $\lambda_1\geq \ldots\geq \lambda_t\geq 0$, $|\lambda| = t$, where
$|\lambda|:=\sum_{i=1}^t\lambda_i$, and the length
$\ell(\lambda)$ is the number of nonzero parts of $\lambda$. Note that we add
and suppress zero entries in $\lambda$ without further notice, so that we can
also compare partitions of different lengths. For partitions $\lambda,\lambda'$
of integers $t,t'$, respectively, we denote $\lambda \le \lambda'$ if and only
if $\lambda_i\leq \lambda_i'$, for all $i=1,\ldots,\ell(\lambda)$.

The space of polynomials of degree at most $t$ on $\G_{k,d}$ is given by
\begin{equation*}
  \Pol_{t}(\G_{k,d}):=   \{ f|_{\G_{k,d}}: f\in\C[X]_t\},
\end{equation*}
where $\C[X]_t$ is the set of polynomials of degree at most
$t$ in $d^{2}$ many variables arranged as a matrix $X\in\C^{d\times d}$, and
$f|_{\G_{k,d}}$ denotes the restriction of $f$ to $\G_{k,d}$. This polynomial
space decomposes into
\begin{equation*}
  \Pol_{t}(\G_{k,d}) = \!\!\!\!\!\!\bigoplus_{\begin{smallmatrix}
      |\lambda| \le t,\\
      \ell(\lambda) \le \min\{k,d-k\}
    \end{smallmatrix}} \!\!\!\!\!\!H_{\lambda}(\G_{k,d}),
\end{equation*}
so that its dimension is calculated by adding the dimensions of each of the occurring $\mathcal H_{2\lambda}^{d}$. The dimension of $\mathcal H_{2\lambda}^{d}$ is specified in \cite[Formulas
  (24.29) and (24.41)]{Fulton:1991fk}. 
%


\section{Polynomials on unions of Grassmannians}\label{sec:uni}

Given a non-empty set $\mathcal{I}\subset\{1,\ldots,d-1\}$, the 
corresponding union of Grassmannians is
\begin{equation*}
  \mathcal{G}_{\mathcal{I},d} :=\bigcup_{k\in\mathcal{I}} \mathcal{G}_{k,d} = \{ P \in \R^{d\times d}_{\sym} \;:\; P^{2}=P ,\; \Tr(P) \in \mathcal{I}\},\qquad \G_{d}:=\bigcup_{k=1}^{d-1} \G_{k,d}.
\end{equation*}
An orthogonally invariant measure on $\mathcal{G}_{\mathcal{I},d}$ is derived by the sum of the corresponding measures on the single Grassmannians. In Section \ref{sec:cub}, we shall also allow for weighted sums. According to \eqref{eq:decomp L2 single}, the corresponding space of complex-valued, square-integrable
functions $L^{2}(\G_{\mathcal{I},d})$ decomposes into
\begin{equation}\label{eq:L2 union abstract}
  L^{2}(\G_{\mathcal{I},d}) = \bigoplus_{\lambda\in\Lambda^d_\mathcal{I}} H_{\lambda}(\G_{\mathcal{I},d}), \qquad H_{\lambda}(\G_{\mathcal{I},d}) \cong |\mathcal{I}^d_\lambda |\mathcal H_{2\lambda}^{d},
\end{equation}
where $\Lambda^d_\mathcal{I}$ is the set of all partitions $\lambda$ of length at most $\max_{k \in \mathcal{I}}( \min(k, d-k))$ and the multiplicities are the cardinality of  
\begin{equation*}
  \mathcal{I}^d_\lambda := \{ k \in \mathcal{I} : \ell(\lambda) \le \min(k,d-k) \}.
\end{equation*}
As for a single Grassmannian, we consider polynomials on $\G_{\mathcal{I},d}$ given by
multivariate polynomials in the matrix entries of a given projector
$P \in \G_{\mathcal{I},d}$, i.e., 
\begin{equation*}
  \Pol_{t}(\G_{\mathcal{I},d}) := \{ f|_{\G_{\mathcal{I},d}} : f\in\C[X]_t \}.
\end{equation*}
This space decomposes orthogonally into
  \begin{equation*}
    \Pol_{t}(\G_{\mathcal{I},d}) =\bigoplus_{|\lambda|\le t,\;\lambda\in\Lambda^d_\mathcal{I}} H_{\lambda}^{t}(\G_{\mathcal{I},d}), \qquad H_{\lambda}^{t}(\G_{\mathcal{I},d}) \cong   \mu_{\lambda}^{d}(\mathcal{I},t) \mathcal H_{2\lambda}^{d},
  \end{equation*}
 where the multiplicities $ \mu_{\lambda}^{d}(\mathcal{I},t)$ still need to be determined. Indeed, this is the topic of the first part of the present paper. 
\begin{theorem}\label{the:erstes d}
For $t \ge 0$ and $\mathcal{I} = \{ k_i\}_{i=1}^r \subset \{1,\dots,d-1\}$ with $r=|\mathcal{I}|$, the multiplicity of $\mathcal{H}^d_{2\lambda}$ in the direct sum $\bigoplus_{i=1}^{s}
    \Pol_{t-i+1}(\G_{k_{i},d})$ is a lower bound for that in $ \Pol_{t}(\G_{\mathcal{I},d})$, where $s:=\min\{t+1,|\mathcal{I}|\}$.
\end{theorem}
 \begin{proof}
The cases $t=0$ and $|\mathcal{I}|=1$ are trivially fulfilled. Suppose $t\geq 1$ and $r\geq 2$. The restriction mapping $|_{\mathcal{G}_{k_1,d}}: \Pol_{t}(\G_{\mathcal{I},d})\rightarrow\Pol_{t}(\G_{k_1,d})$ is orthogonally invariant and surjective, so that $\Pol_{t}(\G_{\mathcal{I},d})$ is equivalent to $\Pol_{t}(\G_{k_1,d}) \oplus \Null(|_{\mathcal{G}_{k_1,d}})$. Let $\Tr(\cdot)$ denote the trace. Since $(\Tr(\cdot) - k_1) \Pol_{t-1}(\G_{\mathcal{I},d}) $ is equivalent to $\Pol_{t-1}(\G_{\{k_2,\ldots,k_r\},d})$ and
\begin{equation*}
(\Tr(\cdot) - k_1) \Pol_{t-1}(\G_{\mathcal{I},d}) \subset \Null(|_{\mathcal{G}_{k_1,d}}),
\end{equation*}
the number of irreducible components $\mathcal{H}^{d}_{2\lambda}$ in $\Pol_{t}(\G_{\mathcal{I},d})$ is bigger or equals the respective number in $\Pol_{t}(\G_{k_1,d})\oplus \Pol_{t-1}(\G_{\{k_2,\ldots,k_r\},d})$. An induction over $t$ and $r$ completes the proof.
\end{proof}
If we order $\mathcal{I} =\{ k_i\}_{i=1}^r$ by $
    \min\{k_{1},d-k_{1}\} \ge \dots \ge \min\{k_{r},d-k_{r}\}$, 
  then counting the actual occurrences of $\mathcal{H}_{2\lambda}^d$ in $\bigoplus_{i=1}^{s}
    \Pol_{t-i+1}(\G_{k_{i},d})$ yields the following explicit lower bound.
\begin{corollary}\label{cor:11}
For $t \ge 0$ and $\mathcal{I} \subset \{1,\dots,d-1\}$, the multiplicities $ \mu_{\lambda}^{d}(\mathcal{I},t)$ satisfy
\begin{equation}\label{eq:est lower}
    \mu^d_{\lambda}(\mathcal{I},t) \geq  
      \min\{ t - |\lambda|+1,|\mathcal{I}^d_\lambda|\}, \quad  0\le|\lambda|\le t,  \quad \lambda\in\Lambda^d_\mathcal{I}.
  \end{equation}
\end{corollary}
Reformulation yields $\mu^d_{\lambda}(\mathcal{I},|\lambda|+s) \geq  
      \min\{ s+1,|\mathcal{I}^d_\lambda|\}$, for $s\geq 0$. 
Due to \eqref{eq:L2 union abstract}, the upper bound $\mu_{\lambda}^{d}(\mathcal{I},t)\leq |\mathcal{I}^d_\lambda|$ holds. Subsequent sections shall reveal that equality holds in \eqref{eq:est lower}. However, this requires a closer look at relations among irreducible representations and their reproducing kernels, see Appendix \ref{app:rep k} for some basics on reproducing kernels that shall be used in the following.

\section{Determining the multiplicities for few special cases
}\label{sec:H}
The space of polynomials on
$\R^{d\times d}_{\sym}$ of degree at most $t$ 
and its subspace of homogeneous polynomials of degree
$t$ are denoted by 
\begin{align*}
\Pol_t(\R^{d\times d}_{\sym}) &: = \{ f|_{\R^{d\times d}_{\sym}} : f\in \C[X]_t\},\\
\Hom_t(\R^{d\times d}_{\sym})&:=\{f\in \Pol_t(\R^{d\times d}_{\sym}) : f(\alpha X)=\alpha^t f(X), \; \alpha\in\R,\;X\in \R^{d\times d}_{\sym}\},
\end{align*}
respectively. The differential inner product between $f,g \in \Pol_{t}(\R^{d\times
  d}_{\sym})$ 
given by 
\[
(f,\overline g)_{\mathrm D} := f(D) \overline g(0),\qquad  \text{where }
\mathrm D:= \Big( \frac12 ( 1 + \delta_{i,j} ) \partial_{i,j} \Big)_{i,j=1,\ldots,d}\;,
\]
is orthogonally invariant, cf.~\cite{That:1976mz}, inducing the orthogonal
decomposition
\begin{equation*}
\Pol_t(\R^{d\times d}_{\sym})  = \bigoplus_{s=0}^t \Hom_s(\R^{d\times d}_{\sym}).
\end{equation*}

%
%
\begin{remark}
\label{rem:trace}
In accordance with \cite{Gross:1987bf}, the mapping $(X,Y)\mapsto \frac{1}{s!}\Tr(XY)^s$ 
is the reproducing kernel for $\Hom_s(\R^{d\times d}_{\sym})$ with respect to
the differentiation inner product, cf.~also Appendix \ref{app:rep k}. 
\end{remark}
Let $\GL(d)$ denote the general linear group. The space $\Hom_t(\R^{d\times d}_{\sym}) $ decomposes orthogonally into subspaces $F_\lambda(\R^{d\times d}_{\sym})$ invariant under the action $f \mapsto f(L \cdot L^{\top})$, for $L\in\GL(d)$, by
\begin{equation}\label{eq:hom in F}
\Hom_t(\R^{d\times d}_{\sym}) = \bigoplus_{|\lambda|=t,\;\;\ell(\lambda)\leq d} F_\lambda(\R^{d\times d}_{\sym}),\qquad F_\lambda(\R^{d\times d}_{\sym})\cong \mathcal F_{2\lambda}^{d},
\end{equation}
where $\mathcal F_{2\lambda}^{d}$ is the irreducible
representation of $\GL(d)$ associated to $2\lambda=(2\lambda_{1},\dots,2\lambda_{t})$, cf.~\cite{Gross:1987bf}.

 \begin{remark}
One can check that the function $(X,Y)\rightarrow \frac{1}{|\lambda|!} C_\lambda(X Y)$ 
is the reproducing kernel for $F_\lambda(\R^{d\times d}_{\sym})$ with respect to the differentiation inner product, where $C_\lambda$ is the zonal polynomial of index $\lambda$, cf.~\cite{Gross:1987bf} and Appendix \ref{sec:zon}.
\end{remark}
By restricting the group action from $\GL(d)$ to $\mathcal{O}(d)$, the space $F_\lambda(\R^{d\times d}_{\sym})$ decomposes
further as
\begin{equation}\label{eq:F in Hs}
  F_{\lambda}(\R^{d\times d}_{\sym}) \cong \mathcal F_{2\lambda}^{d} \cong \bigoplus_{\lambda'\leq 2\lambda}\nu^d_{2\lambda,\lambda'} \mathcal{H}^d_{\lambda'},\qquad \ell(\lambda) \le d,
\end{equation}
where the multiplicities
$\nu^d_{2\lambda,\lambda'}\in\N_0$ are determined by the corresponding branching rule,
cf.~\cite{Howe:2004jt, Koike:1987ez, Fulmek:ss}. If $\lambda'$ is such that $\mathcal{H}^d_{\lambda'}$ is not defined, then we simply put $\nu^d_{2\lambda,\lambda'}=0$. 
%
%
One observes 
\begin{equation*}
\nu^d_{2\lambda,(0)}=1,\quad\text{for $\ell(\lambda)\le d$},\qquad\qquad \nu^d_{2\lambda,2\lambda}=1,\quad\text{for
$\ell(\lambda) \le d/2$.}
\end{equation*}
Thus, we obtain $\mu_{\lambda}^{d}(\mathcal{I},t)  = 1$, for $|\lambda|=t\geq 0$ with $\lambda\in\Lambda^d_\mathcal{I}$. The latter enables us to verify that the lower bounds in Corollary \ref{cor:11} are an equality for $|\mathcal{I}|=2$. Indeed, if we order $\mathcal I = \{ k_{1}, k_{2} \} \subset \{1,\dots,d-1\}$ by $\min\{k_{1},d-k_{1}\}  \ge \min\{k_{2},d-k_{2}\}$, then Theorem~\ref{the:erstes d} and $\mu_{\lambda}^{d}(\mathcal{I},t)  = 1$ imply
\begin{equation*}
 \mu_{\lambda}^{d}(\{k_1,k_2\},t) = \begin{cases}  1 , & |\lambda|=t \text{ with } \ell(\lambda)\leq \min\{k_{1},d-k_{1}\},\\
  1 , & |\lambda|\leq t-1 \text{ with } \min\{k_{2},d-k_{2}\}<\ell(\lambda)\leq \min\{k_{1},d-k_{1}\},\\
 2,& |\lambda|\leq t-1 \text{ with } \ell(\lambda)\leq \min\{k_{2},d-k_{2}\},
 \end{cases}
\end{equation*}
which means that equality holds in \eqref{eq:est lower}. 

For arbitrary $\mathcal{I} \subset \{1,\dots,d-1\}$, we can still use \eqref{eq:F in Hs} to determine the multiplicities $\mu_\lambda^d(\mathcal{I},t)$ provided that $t$ is sufficiently small, but this takes some preparation. The space of homogeneous polynomials restricted to
$\G_{\mathcal{I},d}$ is denoted by
$ 
\Hom_t(\mathcal{G}_{\mathcal{I},d})  := \{ f|_{\mathcal{G}_{\mathcal{I},d}} : f\in \Hom_t(\R^{d\times d}_{\sym}) \}.
$ 
It is important to notice that the restriction of $\Hom_t(\R^{d\times d}_{\sym})$ to $\mathcal{G}_{\mathcal{I},d}$ yields (almost) the entire space $\Pol_t(\mathcal{G}_{\mathcal{I},d})$:
\begin{theorem}\label{the:structure} For $t\ge 0$ and $\mathcal{I} \subset \{1,\dots,d-1\}$, the polynomial space decomposes into
\begin{align*}
\Pol_{t}(\mathcal{G}_{\mathcal{I},d}) & =
\begin{cases}
  \Hom_t(\mathcal{G}_{\mathcal{I},d})\oplus \Hom_0(\G_{\mathcal{I},d}), &1\le t \leq |\mathcal I|-1,\\
  \Hom_t(\mathcal{G}_{\mathcal{I},d}), &\text{else}.
\end{cases} 
\end{align*}
\end{theorem}
\begin{proof}
  First we note that $\Tr(X \cdot)^{t}$, $X \in \R^{d\times d}_{\sym}$, linearly generates the space
  $\Hom_{t}(\G_{\mathcal{I},d})$, cf.~Remark~\ref{rem:trace}, \eqref{eq:S0 def}, and \eqref{eq:HKYspann}. Now, for $t\ge 1$ and
  $\mathcal{I} \subset \{1,\dots,d-1\}$, we observe
 \[
  \Tr(XP)^{s}=\Tr(XP)^{s-1}\Tr(XP^{t-s+1}),
  \qquad X \in \R^{d\times d}_{\sym},\quad P \in \G_{\mathcal{I},d}, \quad 1 \le
  s \le t.
  \]
Since the term on the right hand side is a homogeneous polynomial of degree
  $t$ in $P$ restricted to $\G_{\mathcal{I},d}$, we deduce
  \begin{equation}\label{eq:hom et}
  \Hom_{t}(\G_{\mathcal{I},d}) = \{ f|_{\mathcal{G}_{\mathcal{I},d}} : f\in
  \Hom_s(\R^{d\times d}_{\sym}),\; s=1,\ldots,t \}.
 \end{equation}
Thus, it remains to check that $1|_{\G_{\mathcal{I},d}} \not\in \Hom_{t}(\G_{\mathcal{I},d})$ if and only if $1\le t \le |\mathcal I|- 1$. 
  Since $1|_{\G_{\mathcal{I},d}}$ is orthogonally
  invariant, it is sufficient to consider the orthogonally invariant subspace of
  $\Hom_{t}(\G_{\mathcal{I},d})$ denoted by $\Hom_{t}^{\mathcal O(d)}(\G_{\mathcal{I},d})$. 
 For $f \in \Hom_{t}(\R^{d\times d}_{\sym})$, we define 
  \begin{equation*}
  f_{\circ}(X) := \int_{\mathcal O(d)} f(O X O^{\top})d\mu_{\mathcal
    O(d)}(O),\qquad X \in \R^{d\times d}_{\sym},
  \end{equation*}
so that $\Hom_{t}^{\mathcal O(d)}(\G_{\mathcal{I},d}) = \spann \big\{ f_{\circ}|_{\mathcal G_{\mathcal{I},d}} : f \in
    \Hom_{t}(\R^{d\times d}_{\sym}) \big\}$. 
According to invariant theory, the
  ring of orthogonally invariant polynomials on $\R^{d\times d}_{\sym}$ is generated by
  polynomials of the form $\Tr(X^{l})$, $X \in \R^{d\times d}_{\sym}$, ${l} \in \N_{0}$, cf.~\cite[Theorem~7.1]{Procesi:1976ul}. 
  Since $f_{\circ}$ is also homogeneous of degree $t$, we observe
    \begin{equation*}
  f_{\circ}(X)= \sum_{
      {l}_{1}+\dots +{l}_{t}=t
  }
  f_{{l}_{1},\dots,{l}_{t}}\prod_{i=1}^{t} \Tr(X^{{l}_{i}}),\quad X\in \R^{d\times d}_{\sym}.
  \end{equation*}
For $t\geq 1$, the restriction $f_{\circ}|_{\G_{\mathcal{I},d}}$ is a linear combination of the functions $\Tr(\cdot)^{s}|_{\G_{\mathcal{I},d}}$, $s=1,\dots,t$. According to \eqref{eq:hom et}, these functions are contained in $\Hom_{t}(\G_{\mathcal{I},d})$,
  so that $\Hom_{t}^{\mathcal O(d)}(\G_{\mathcal{I},d})$ is spanned by $\Tr(\cdot)|_{\mathcal{G}_{\mathcal{I},d}}$, \ldots, $\Tr(\cdot)^t|_{\mathcal{G}_{\mathcal{I},d}}$, for $t\ge 1$. 
The invertibility of the Vandermonde matrix implies that $1_{\mathcal{G}_{\mathcal{I},d}}$ and
  $\Tr(\cdot)^{i}|_{\mathcal{G}_{\mathcal{I},d}}$, $i=1,\dots,t$, are linearly
  independent if and only if $1\leq t \le |\mathcal I|- 1$. 
\end{proof}

\begin{remark}\label{rem:erster}
Fix some $t\geq 0$ and $\mathcal{I} \subset \{1,\dots,d-1\}$. The map
$ 
 (P,Q)\mapsto \Tr(PQ)^t
$ 
is a reproducing kernel for $\Hom_{t}(\G_{\mathcal{I},d})$. For any constant $C > 0$, the map $ 
 (P,Q)\mapsto \Tr(PQ)^t+C
$ 
is a reproducing kernel for $\Pol_{t}(\G_{\mathcal{I},d})$.
\end{remark}
Next, we determine the multiplicities $\mu_\lambda^d(\mathcal{I},t)$, for $t=1,2,3$, by deriving upper bounds that match the lower bounds in Corollary \ref{cor:11}. The decompositions \eqref{eq:hom in F} and \eqref{eq:F in Hs} yield $\Hom_{1}(\R^{d\times d}_{\sym}) \cong \mathcal H_{(0)}^{d} \oplus \mathcal H_{(2)}^{d}$, for $d\ge 2$, and 
  \begin{equation}\label{eq:branch rules explit}
    \Hom_{2}(\R^{d\times d}_{\sym}) \cong \begin{cases} 
      2\mathcal H_{(0)}^{d} \oplus 2\mathcal H_{(2)}^{d} \oplus \mathcal H_{(4)}^{d} \oplus \mathcal H_{(2,2)}^{d}, & d\ge 4,\\
      2\mathcal H_{(0)}^{d} \oplus 2\mathcal H_{(2)}^{d} \oplus \mathcal H_{(4)}^{d}, & d = 3,\\
      2\mathcal H_{(0)}^{d} \oplus \;\;\mathcal H_{(2)}^{d} \oplus \mathcal H_{(4)}^{d}, & d = 2.\\
    \end{cases}
    \end{equation}
The multiplicities of $\mathcal{H}_\lambda^d$ in $\Hom_{t}(\R^{d\times d}_{\sym})$ are upper bounds for $\mu^d_{\lambda}(\mathcal{I},t)$ since the restriction mapping is orthogonally invariant. For $t=1,2$, the lower
  bounds in Corollary \ref{cor:11} are matched, so that we have determined $\mu^d_{\lambda}(\mathcal{I},1)$ and $\mu^d_{\lambda}(\mathcal{I},2)$ for any index set
  $\mathcal{I} \subset \{1,\dots,d-1\}$. For $t=3$, the analysis is more difficult, and we observe that the branching rules yield
    \begin{align*}
 \Hom_{3}(\R^{d\times d}_{\sym}) \cong
    \begin{cases} 3\mathcal H_{(0)}^{d} \oplus 4\mathcal H_{(2)}^{d} \oplus 2\mathcal H_{(4)}^{d} \oplus 2\mathcal H_{(2,2)}^{d} \oplus \mathcal H_{(3,1)}^{d} \oplus \mathcal H_{(6)}^{d} \oplus \mathcal H_{(4,2)}^{d} \oplus \mathcal H_{(2,2,2)}^{d}, & d \ge 6,\\
 3\mathcal H_{(0)}^{d} \oplus 4 \mathcal H_{(2)}^{d} \oplus 2\mathcal H_{(4)}^{d} \oplus 2\mathcal H_{(2,2)}^{d} \oplus \mathcal H_{(3,1)}^{d} \oplus \mathcal H_{(6)}^{d} \oplus \mathcal H_{(4,2)}^{d}, & d = 5,\\
 3\mathcal H_{(0)}^{d} \oplus 4 \mathcal H_{(2)}^{d} \oplus 2\mathcal H_{(4)}^{d} \oplus \;\;\mathcal H_{(2,2)}^{d} \oplus \mathcal H_{(3,1)}^{d} \oplus \mathcal H_{(6)}^{d} \oplus \mathcal H_{(4,2)}^{d}, & d = 4,\\
 3\mathcal H_{(0)}^{d} \oplus 3 \mathcal H_{(2)}^{d} \oplus 2\mathcal H_{(4)}^{d} \qquad\qquad\,\,\oplus \mathcal H_{(3,1)}^{d} \oplus \mathcal H_{(6)}^{d}, & d = 3,\\
 2\mathcal H_{(0)}^{d} \oplus 2\mathcal H_{(2)}^{d} \oplus \;\;\mathcal H_{(4)}^{d} \qquad\qquad\qquad\qquad\,\oplus \mathcal H_{(6)}^{d}, & d = 2.
    \end{cases}
  \end{align*}
  The multiplicity of $\mathcal H_{(2)}^{d}$ in $\Hom_{3}(\R^{d\times d}_{\sym}) $ does not match the lower bound in Corollary \ref{cor:11}. 
Instead, we found the kernel
\begin{align}\label{eq:KKK}
\begin{split}
  K(X,Y) = & \frac{1}{d+2} \Big(\Tr(X^{2}Y^{2}) \Tr(X Y) - \Tr(X^{2}Y)\Tr(XY^{2})\Big)\\
  & - \frac{1}{(3d+4)(d+2)} \Big(\Tr(X^{2}Y^{2})\Tr(X)\Tr(Y) - \Tr(X^{2}Y)\Tr(X)\Tr(Y^{2}) \\
  & \hspace{2.65cm} - \Tr(X Y^{2})\Tr(X^{2})\Tr(Y) + \Tr(XY)\Tr(X^{2})\Tr(Y^{2})
  \Big),
  \end{split}
\end{align}
which may not have been observed in the literature yet, reproduces a subspace of $\Hom_{3}(\R^{d\times d}_{\sym}) $ equivalent to $\mathcal{H}^d_{(3,1)}\oplus \mathcal{H}^d_{(2)}$, for $d>2$, and equivalent to 
$\mathcal H_{(2)}^{d}$, for $d=2$, respectively. 
Since $K$ vanishes on any
Grassmannian, i.e., $K(X,Y)=0$, for all $Y\in \G_{k,d}$, $k \in \{1,\dots,d-1\}$
and $X\in\R^{d\times d}_{\sym}$, we deduce that the multiplicity of $\mathcal H_{(2)}^{d}$ in
$\Pol_{3}(\G_{\mathcal{I},d})$ is less than in
$ \Hom_{3}(\R^{d\times d}_{\sym}) $. Now, for even partitions, the resulting upper bounds on the multiplicities match the lower bounds
in Corollary \ref{cor:11} for $t=3$ and any index set
$\mathcal{I} \subset \{1,\dots,d-1\}$.

To determine the multiplicities of the irreducible subspaces of $\Pol_{3}(\G_{\mathcal{I},d})$ for $|\mathcal{I}|\geq 3$, we used the reproducing kernel $K$ in \eqref{eq:KKK}. We shall further explore the reproducing kernels of the irreducible components to determine the multiplicities in $\Pol_t(\mathcal{G}_{\mathcal{I},d})$ for $|\mathcal{I}|\geq 3$ and $t\geq 4$.

\section{Zonal kernels and harmonic analysis on unions of Grassmannians}\label{sec:L2kernels} 
For $\ell(\lambda)\leq \min(k,d-k,l,d-l)$, the spaces $H_{\lambda}(\G_{k,d})$ and $H_{\lambda}(\G_{l,d})$ are equivalent, hence, there is a real intertwining isomorphism $T_{\lambda}^{k,l}: H_{\lambda}(\G_{k,d}) \to H_{\lambda}(\G_{l,d})$. In particular, $T_{\lambda}^{k,l}$ commutes with complex conjugation and the group action. It can be realized by an integral transform with a unique real-valued zonal function $p_{\lambda}^{k,{l}}:\G_{k,d}\times\G_{{l},d}\to \R$, so that 
\begin{equation}\label{eq:int func}
T_{\lambda}^{k,{l}}f  = \int_{\mathcal{G}_{{k},d}} p_{\lambda}^{k,{l}}(P,\cdot)f(P)\d\sigma_{k,d}(P).
\end{equation} 
Note that zonal means $p_{\lambda}^{k,{l}}(OPO^\top,OQO^\top)=p_{\lambda}^{k,{l}}(P,Q)$, for $O\in\mathcal{O}(d)$ and $P\in\mathcal{G}_{{k},d}$, $Q\in\mathcal{G}_{{l},d}$. 
For fixed  $1 \le k \le {l} \le \tfrac d 2$ and $\ell(\lambda)\leq k$, the intertwining functions $ p_{\lambda}^{k,{l}}$ were studied in \cite{James:1974aa} and expanded into the zonal polynomials by
\begin{equation}
    \label{eq:pkl_jacobi}
    p_{\lambda}^{k,{l}}(P,Q) = b_{\lambda}^{k,{l},d} \sum_{\lambda'\le \lambda}
  c_{\lambda,\lambda'}^{\frac d 2} q_{\lambda,\lambda'}(\tfrac{k}{2}) q_{\lambda,\lambda'}(\tfrac {l} 2)
  C_{\lambda'}(PQ),\qquad \text{ for $\ell(\lambda)\leq k$, \;\;$1 \le k \le {l} \le \tfrac d 2$,}
   \end{equation}
where $b_{\lambda}^{k,{l},d}\in\R$ is a scaling constant, $c_{\lambda,\lambda}^{\frac d 2}=1$, and $q_{\lambda,\lambda'}$ is a polynomial of degree $|\lambda|-|\lambda'|$ given by 
\begin{equation*}
  q_{\lambda,\lambda'}(x):= \prod_{i=1}^{m} (x-\tfrac 1
  2(i-1)+\lambda'_{i})_{\lambda_{i}-\lambda'_{i}} = \frac{(x)_{\lambda}}{(x)_{\lambda'}}, \quad x
  \in \R.
\end{equation*}
Potential zeros in the denominator of $\frac{(x)_{\lambda}}{(x)_{\lambda'}}$ cancel
out, so that the fraction is well-defined. Up to the scaling, which we have not specified yet, the functions $p_{\lambda}^{k,k}$ are the reproducing kernels for $H_{\lambda}(\G_{k,d})$ with respect to the $L^2$ inner product when $1\leq k\leq d/2$ and $\ell(\lambda)\leq k$, cf.~\cite{James:1974aa}. 

The sum of the right hand side in \eqref{eq:pkl_jacobi} is still well-defined for all $k,l\in\{1,\ldots,d-1\}$ and $\ell(\lambda)\leq \min(k,d-k,l,d-l)$. One of our contributions going beyond \cite{James:1974aa} is to determine the reproducing kernel of $H_{\lambda}^{|\lambda|}(\G_{d})$ with the help of a particular extension of \eqref{eq:pkl_jacobi} to this broader range of parameters. Recall that $\G_{d}=\bigcup_{k=1}^{d-1} \G_{k,d}$.
\begin{theorem}\label{the:first and foremost}
Let $\ell(\lambda)\leq d/2$. The reproducing kernel of $H_{\lambda}^{|\lambda|}(\G_{d})$ with respect to the $L^2$ inner product is a multiple of 
\begin{equation*}
 p_{\lambda}(P,Q) := 
    \sum\limits_{\lambda' \le \lambda} c_{\lambda,\lambda'}^{\frac d 2} q_{\lambda,\lambda'}(\tfrac
    12 \Tr(P)) q_{\lambda,\lambda'}(\tfrac 12 \Tr(Q)) C_{\lambda'}(P
    Q),\quad P, Q \in \G_{d}.
   \end{equation*}
\end{theorem}
In order to verify Theorem~\ref{the:first and foremost}, we choose a suitable normalization of the intertwining functions, induced by the following selection of intertwining operators.
\begin{proposition}\label{prop:11}
For all $k,l\in \{1,\ldots,d-1\}$ and $\ell(\lambda)\leq \min(k,d-k,l,d-l)$, there are real isometric isomorph 
 intertwining operators $T_{\lambda}^{k,{l}}:H_{\lambda}(\G_{{k},d}) \to H_{\lambda}(\G_{l,d})$ such that the following diagram commutes, 
\begin{center}
\begin{minipage}[c]{.45\textwidth}
\begin{tikzpicture}
  \matrix (m) [matrix of math nodes,row sep=3em,column sep=8em,minimum width=2em]
  {
    H_{\lambda}(\G_{k,d}) & \\
     H_{\lambda}(\G_{m,d}) & H_{\lambda}(\G_{l,d}),\\};
  \path[-stealth]
    (m-1-1) edge node [left] {$T_{\lambda}^{k,m}$ } (m-2-1)
	 edge node [right] {\hspace{3ex}$T_{\lambda}^{k,l}$} (m-2-2)
             (m-2-1)   edge  node [above] {\hspace{-5ex}$T_{\lambda}^{m,l}$} (m-2-2);
\end{tikzpicture}
\end{minipage}
for all $k,m,l\in \{\ell(\lambda),\dots,d-\ell(\lambda)\}$.
\end{center}
\end{proposition}
\begin{proof}
Let us fix an index $s$ with $\ell(\lambda)\leq s\leq d - \ell(\lambda)$. There are real isometric isomorph operators 
  $\tilde T_{\lambda}^{k,s}$, $k \in \{\ell(\lambda),\dots, \ell(\lambda)\}$, which
  intertwine the spaces $H_{\lambda}(\G_{k,d})$ and $H_{\lambda}(\G_{s,d})$. We now define
  \begin{equation*}
    T_{\lambda}^{k,{l}} := (\tilde{T}_{\lambda}^{{l},s})^*\tilde{T}_{\lambda}^{k,s},\qquad \ell(\lambda)\leq k,{l}\leq d-\ell(\lambda),
  \end{equation*}
and straightforward calculations yield the statement.   
\end{proof}
The integral operators in Proposition~\ref{prop:11} induce intertwining functions $p_\lambda^{k,{l}}$ via \eqref{eq:int func} satisfying 
  \begin{align}
 \dim(\mathcal H_{2\lambda}^{d}) & =\int_{\G_{k,d}}\int_{\G_{{l},d}} |p_{\lambda}^{k,{l}}(P,Q)|^{2} \d\sigma_{{l},d}(Q) \d\sigma_{k,d}(P),     \label{eq:normalizedplk} \\
p_{\lambda}^{k,{l}}(P,Q) & =    \big(p_{\lambda}^{k,m}(P,\cdot), p_{\lambda}^{m,{l}}(\cdot,Q)\big)_{\G_{m,d}} 
, \quad\text{for $m \in \{\ell(\lambda),\dots,d-\ell(\lambda)\}$. }\label{eq:pkl_rel}
  \end{align}
\begin{remark}
The intertwining functions $p_{\lambda}^{k,k}$ are the reproducing kernels of $H_\lambda(\mathcal{G}_{k,d})\cong \mathcal{H}^d_{2\lambda}$ with respect to the standard $L^2$ inner product for $k\in\{1,\ldots,d-1\}$ and $\ell(\lambda)\leq \min(k,d-k)$.  
\end{remark}

Let $P_{\lambda}^{k,{l}}:\G_{\mathcal{I},d}\times \G_{\mathcal{I},d} \to \R$ denote the zero extension of $p_\lambda^{k,{l}}$. It follows from \cite{James:1974aa} that the collection of zonal functions $\{P_{\lambda}^{k,{l}} : \lambda\in\Lambda^d_\mathcal{I};\; k,{l}\in\mathcal{I}^d_\lambda\}$ is an orthogonal basis for $L^{2}_{\mathcal O(d)}(\G_{\mathcal{I},d}\times \G_{\mathcal{I},d})$, the space of square integrable functions that are zonal. Hence, any zonal function $f \in L^{2}_{\mathcal O(d)}(\G_{\mathcal{I},d} \times \mathcal G_{\mathcal{I},d})$ can be
expanded into a Fourier series, i.e.,
\begin{equation}\label{FK exp a}
f = \sum_{\lambda\in\Lambda^d_\mathcal{I} } \sum_{
    k,{l} \in \mathcal{I}^d_\lambda
} \hat f_{\lambda}^{k,{l}} P_{\lambda}^{k,{l}} = \sum_{\lambda\in\Lambda^d_\mathcal{I}} \Tr( \hat
f_{\lambda}^{\top} P_{\lambda}),
\end{equation}
where the Fourier coefficients $\hat f_{\lambda}^{k,{l}}$ and the basis functions $P_{\lambda}^{k,{l}}$ are arranged in matrix
form
\begin{equation*}
  \hat f_{\lambda}  := \big(\hat f_{\lambda}^{k,{l}} \big)_{k,{l} \in \mathcal{I}^d_\lambda} \in \C^{|\mathcal{I}^d_\lambda|\times |\mathcal{I}^d_\lambda|}, \qquad
  P_{\lambda}  := \big( P_{\lambda}^{k,{l}} \big)_{k,{l} \in \mathcal{I}^d_\lambda}.
\end{equation*}
Convolving two continuous zonal functions
$f,g \in L^{2}_{\mathcal O(d)}(\G_{\mathcal{I},d} \times \G_{\mathcal{I},d})$,
\[
(f * g) (P,Q) := \big(f (P,\cdot), \overline{g (\cdot,Q)}\big)_{\G_{\mathcal{I},d}}= \sum_{k \in \mathcal I}\int_{\G_{k,d}}  \!\!\!\!\!\!f (P,R)g (R,Q)d\sigma_{k,d}(R),
\] 
yields again a continuous zonal function $f*g$. It is straight-forward to check that its Fourier coefficients are $\widehat{(f*g)}_\lambda = \hat f_{\lambda} \hat g_{\lambda} \in \C^{|\mathcal{I}^d_\lambda|\times
 |\mathcal{I}^d_\lambda|}$, $\lambda\in\Lambda^d_\mathcal{I}$. 
 \begin{remark}
This convolution property implies that the kernel $(P,Q)\mapsto \Tr(P_{\lambda}(P,Q))$, for $\lambda\in\Lambda^d_\mathcal{I}$, is the reproducing kernel of $H_{\lambda}(\G_{\mathcal{I},d})$ with respect to the $L^{2}(\G_{\mathcal{I},d})$ inner product.
\end{remark}
 The Fourier
coefficients $\hat K_{\lambda} \in \C^{|\mathcal{I}^d_\lambda|\times |\mathcal{I}^d_\lambda|}$ of a positive definite zonal kernel
$K \in L^{2}_{\mathcal O(d)}(\G_{\mathcal{I},d} \times \G_{\mathcal{I},d})$
are positive semidefinite matrices and thus allow for a spectral decomposition
  \begin{equation}
    \label{eq:hatKpd}
    \hat K_{\lambda} = \sum_{i=1}^{|\mathcal{I}^d_\lambda|} \alpha^i_{\lambda}
    \hat{K}_{\lambda}^i,\qquad \lambda\in\Lambda^d_\mathcal{I},
  \end{equation}
  where $\alpha^1_{\lambda} \ge \dots \ge \alpha^{|\mathcal{I}^d_\lambda|}_{\lambda} \ge 0$   and $\hat{K}_{\lambda}^i\in\C^{|\mathcal{I}^d_\lambda|\times|\mathcal{I}^d_\lambda|}$ are orthogonal rank-$1$ projectors corresponding to an eigenbasis of $\hat{K}_\lambda$. The corresponding kernels 
\begin{equation}\label{eq:new kernels marginals}
K^i_\lambda:= \Tr(\hat{K}^{i\top}_\lambda P_{\lambda}),\quad i=1,\ldots,|\mathcal{I}^d_\lambda|,
\end{equation}
are also positive definite. This spectral decomposition of $K$ yields the irreducible decomposition of the underlying reproducing kernel Hilbert space denoted by $\mathcal{S}(K)$, see also \eqref{eq:S0 def}:
\begin{theorem}
  \label{the:kerneldecomp}
  Let $K : \G_{\mathcal{I},d} \times \G_{\mathcal{I},d} \to \C$ be zonal and positive definite. Under the notations \eqref{eq:hatKpd}, \eqref{eq:new kernels marginals}, we obtain $\mathcal H_{2\lambda}^{d} \cong \mathcal{S}(K^i_\lambda) \subset H_{\lambda}(\G_{\mathcal{I},d})$, for each $i=1,\dots,|\mathcal{I}^d_\lambda|$, and
the multiplicity of $\mathcal H_{2\lambda}^{d}$ in $\mathcal{S}(K)$ equals the rank of the Fourier coefficient $\hat{K}_{\lambda}$, implying the  orthogonal decomposition 
 \begin{equation}
    \label{eq:HKdecomp}
    \mathcal{S}(K) = \bigoplus_{\lambda\in\Lambda^d_\mathcal{I}} \bigoplus_{i=1}^{\rank(\hat
      K_{\lambda})} \mathcal{S}(K^i_\lambda).
  \end{equation}
 \end{theorem}
\begin{proof}
Mercer's Theorem implies that the reproducing kernel Hilbert space $\mathcal{S}(K)$ decomposes
  into the pairwise orthogonal eigenspaces with non-zero eigenvalues associated
  to the integral operator $T_{K}:L^{2}(\G_{\mathcal{I},d}) \to L^{2}(\G_{\mathcal{I},d})$ defined by 
  $
  T_{K} f(P) := (f,K(P,\cdot))_{\G_{\mathcal{I},d}}$, $ P \in \G_{\mathcal{I},d}$, $f\in
  L^{2}(\G_{\mathcal{I},d})$. 
  This decomposition corresponds to the eigenspace decomposition
  of $T_{K}$ in the subspace $H_{\lambda}(\G_{\mathcal{I},d})$. More precisely, 
  the convolution property yields that the kernels $K^i_\lambda$ satisfy, for $i,j=1,\dots,|\mathcal{I}^d_\lambda|$,
\begin{equation*}
  (T_{K} K^i_\lambda(P,\cdot))(Q) = \alpha_{\lambda}^i
  K^i_\lambda(P,Q),\qquad
    (K^i_\lambda(P,\cdot),K^j_\lambda(Q,\cdot))_{\G_{\mathcal{I},d}} =
  \delta_{i,j} K^i_\lambda(P,Q),\qquad P,Q \in \G_{\mathcal{I},d}.
  \end{equation*}
Hence, $K^i_\lambda$ are the reproducing
  kernels for the pairwise orthogonal spaces $\mathcal{S}(K^i_\lambda)$ with respect to the standard inner product.
  The convolution property \eqref{eq:pkl_rel} yields that $\Tr(XP_\lambda(P,\cdot))\in H_{\lambda}(\G_{\mathcal{I},d})$, for any matrix $X\in\C^{|\mathcal{I}^d_\lambda|\times |\mathcal{I}^d_\lambda|}$. Thus, we infer
  $\mathcal{S}(K^i_\lambda) \subset
  H_{\lambda}(\G_{\mathcal{I},d})$, so that 
  \begin{equation*}
  \bigoplus_{i=1}^{|\mathcal{I}^d_\lambda|} \mathcal{S}(K^i_\lambda) \subset
  H_{\lambda}(\G_{\mathcal{I},d}) \cong |\mathcal{I}^d_\lambda| \mathcal H_{2\lambda}^{d},\qquad |\lambda|
  \ge 0.
  \end{equation*}
  Since the spaces $\mathcal{S}(K^i_\lambda) \ne
  \{0\}$,
  $i=1,\dots,|\mathcal{I}^d_\lambda|$, are orthogonally invariant and pairwise
  orthogonal, we obtain $\mathcal{S}(K^i_\lambda) \cong \mathcal
  H_{2\lambda}^d$, which yields \eqref{eq:HKdecomp}.
\end{proof}

Let $K_{2\lambda}:\R^{d\times d}_{\sym} \times \R^{d\times d}_{\sym} \to \R$ denote the reproducing kernel with respect to the differentiation inner product of the irreducible representation $\mathcal H_{2\lambda}^{d}$ in $\Hom_{|\lambda|}(\R^{d\times d}_{\sym})$. One of the ingredients for the following proof of Theorem~\ref{the:first and foremost} is that the restriction $K_{2\lambda}|_{\G_{k,d}\times \G_{{l},d}}$ coincides with $p^{k,l}_\lambda$ up to a multiplicative constant. 
Indeed, we shall follow the strategy in \cite{James:1974aa}:
\begin{proof}[Proof of Theorem~\ref{the:first and foremost}]
 Consider the positive definite zonal kernel
  \[
  K:\G_{d}\times \G_{d} \to \R, \qquad K(P,Q):= C_{\lambda}(P Q),\qquad P,Q
  \in \G_{d}.
  \]
The relation \eqref{eq:HKYspann} implies $\mathcal{S}(K)=F_\lambda(\R^{d\times d}_{\sym})|_{\G_{d}}$. 
 Theorem~\ref{the:kerneldecomp} and  \eqref{eq:F in Hs} yield
  \begin{equation*}
    K = \Tr(\hat K_{\lambda}^{\top} P_{\lambda}) + \sum_{|\lambda'| < |\lambda|,\;\ell(\lambda')\le
      \frac d 2} \Tr(\hat K_{\lambda'}^{\top} P_{\lambda'}).
  \end{equation*}
  Furthermore, \eqref{eq:F in Hs} implies that the irreducible
  representation $\mathcal{H}_{2\lambda}^{d}$ of $\mathcal O(d)$ occurs exactly once in
  $\mathcal{S}(K)$, so that Theorem~\ref{the:kerneldecomp} yields $\mathcal S(v_{\lambda}^{\top} P_{\lambda} v_{\lambda}) \cong H_{\lambda}^{|\lambda|}(\G_{\mathcal{I},d})$, where $\hat K_{\lambda} = v_{\lambda}^{\top}v_{\lambda}$ for some nonzero vector
  $v_{\lambda} \in \R^{|\mathcal{I}^d_\lambda|}$ with $\mathcal{I}=\{1,\ldots,d-1\}$. 
 According to $v_{\lambda}^{\top}P_{\lambda} v_{\lambda}=\Tr(\hat K_{\lambda}^{\top} P_{\lambda})$, we obtain
  \begin{equation}\label{eq:coeff comp 1}
    v_{\lambda}^{\top}P_{\lambda}(P,Q) v_{\lambda} = C_{\lambda}(PQ) - \sum_{|\lambda'| < |\lambda|,\;\ell(\lambda')\le \frac d
      2} \Tr( \hat K_{\lambda'}^{\top} P_{\lambda'}(P,Q)), \qquad P,Q \in \G_{d}.
  \end{equation}
  Note that $v_{\lambda}^{\top}P_{\lambda} v_{\lambda}$ coincides up to a multiplicative factor with the restriction of the kernel $K_{2\lambda}$. We shall verify in the following that the kernel
  $v_{\lambda}^{\top} P_{\lambda} v_{\lambda}$ reflects the expansion of 
  $p_{\lambda}$ into zonal polynomials $C_{\lambda'}$ defined in
  Theorem~\ref{the:first and foremost}. Starting with \eqref{eq:pkl_jacobi}, we exploit the vanishing
  and symmetry properties of the Jacobi polynomials, cf.~Appendix
  \ref{sec:Jacobi}, combined with the symmetry relations of the intertwining
  functions, cf.~Appendix \ref{sec:intertw func}. We observe that, for any
  partition $\lambda$ with $\ell(\lambda) \le \tfrac d 2$ and any
  $(P,Q) \in \G_{k,d}\times G_{l,d}$ with $1\leq k,l\leq d-1$,
  \begin{equation}
    \label{eq:pkl_jacobi_ext}
   p_{\lambda}(P,Q) =   \sum_{\lambda'\le \lambda}
  c_{\lambda,\lambda'}^{\frac d 2} q_{\lambda,\lambda'}(\tfrac{k}{2}) q_{\lambda,\lambda'}(\tfrac {l} 2)
  C_{\lambda'}(PQ)=
    \begin{cases}
      (b_{\lambda}^{k,l,d})^{-1} p_{\lambda}^{k,l}(P,Q), & \ell(\lambda) \le k,l \le
      d-\ell(\lambda),\\
      0, & \text{else},
    \end{cases}
  \end{equation}
  where $b_{\lambda}^{k,l,d} \in \R \setminus \{0\}$. After inserting the expansion from \eqref{eq:pkl_jacobi_ext} into both sides of \eqref{eq:coeff comp 1} via $P_\lambda$ and $P_{\lambda'}$, we aim to compare coefficients of the zonal polynomials. Let $k,l \in \mathcal{I}=\{\ell(\lambda),\dots,d-\ell(\lambda)\}$ be fixed. One can (only) show linear independence of the functions  
  \begin{equation*}
  (P,Q) \mapsto C_{\lambda'}(P Q), \quad (P,Q) \in \G_{k,d} \times \G_{{l},d},\quad \ell(\lambda')\le \min(k,d-k,l,d-l).
  \end{equation*}
Since $\lambda' \le \lambda$ implies $\ell(\lambda') \le \min(k,d-k,l,d-l)$, the zonal polynomials in \eqref{eq:pkl_jacobi_ext} are linearly independent. By applying $P_{\lambda'}(P,Q) = 0$,  $(P,Q) \in \G_{k,d} \times \G_{{l},d}$, for $\ell(\lambda') >\min(k,d-k,l,d-l)$, the summation on the right hand side in \eqref{eq:coeff comp 1} reduces accordingly, and comparing coefficients is justified. Hence, we obtain $1 = v_{\lambda}^{k}v_{\lambda}^{{l}}\, b_{\lambda}^{k,{l},d} c_{\lambda,\lambda}^{\frac d 2}$, so that $c_{\lambda,\lambda}^{\frac d 2}=1$ leads to
\begin{equation*}
p_\lambda|_{\G_{k,d}\times \G_{{l},d}}=v_\lambda^k v_\lambda^l p_\lambda^{k,l},\qquad \ell(\lambda)\leq\min(k,d-k,l,d-l).
\end{equation*}  
In other words, we have verified that, for $\ell(\lambda)\le \Tr(P), \Tr(Q) \le d - \ell(\lambda)$,
  \begin{equation}
    \label{eq:Kpipi_exp}
    v_{\lambda}^{\top}P_{\lambda}(P,Q) v_{\lambda} = \sum_{\lambda' \le \lambda} c^{\frac d 2}_{\lambda,\lambda'}
    q_{\lambda,\lambda'}(\tfrac 12 \Tr(P)) q_{\lambda,\lambda'}(\tfrac 12 \Tr(Q)) C_{\lambda'}(P
    Q)=p_\lambda(P,Q).
  \end{equation}
  For the remaining cases, we observe $v_{\lambda}^{\top}P_{\lambda}(P,Q) v_{\lambda} = 0$, and $p_{\lambda}(P,Q) = 0$ holds due to \eqref{eq:pkl_jacobi_ext}. 
%
%
\end{proof}

For $|\lambda|=0,1,2$, we present the kernels $K_{2\lambda}$ and $p_{\lambda}$ as well as $p^{k,{l}}_{\lambda}$ and $v_{\lambda}^{k}$ in Appendix \ref{app:ex}. We now replace $P_\lambda$ with $p_\lambda$ and still have a suitable Fourier expansion:  
\begin{corollary}
  \label{cor:basischange}
  For $\mathcal{I} \subset \{1,\dots,d -1 \}$, let $K : \G_{\mathcal{I},d} \times \G_{\mathcal{I},d} \to \C$ be a positive definite zonal
  kernel. Then there is a unique symmetric function $\widehat{K}_\lambda:\mathcal{I}^d_\lambda\times\mathcal{I}^d_\lambda\to \C$ such that 
  \[
    K(P,Q) = \sum_{\lambda\in\Lambda^d_\mathcal{I}} \widehat{K}_\lambda(\Tr(P),\Tr(Q)) p_{\lambda}(P,Q),
    \qquad P,Q \in \G_{\mathcal{I},d},
  \]
and the multiplicity of 
  $\mathcal H_{2\lambda}^{d}$ in $\mathcal S(K)$ is the rank of the matrix
  $\big(\widehat{K}_\lambda(k,{l})\big)_{k,{l}\in \mathcal{I}^d_\lambda}$.

\end{corollary}
\begin{proof}
The Fourier expansion \eqref{FK exp a} with coefficients $\hat K_{\lambda}^{k,{l}}$ yields $\widehat{K}_\lambda(k,{l}) = (v_{\lambda}^{k})^{-1}\hat K_{\lambda}^{k,{l}}
    (v_{\lambda}^{{l}})^{-1}$, where $v_{\lambda}$ is given in \eqref{eq:Kpipi_exp}. Therefore, the rank of $\big(\widehat{K}_\lambda(k,{l})\big)_{k,{l}\in \mathcal{I}^d_\lambda}$ is the same as the one of the original Fourier coefficient, so that Theorem~\ref{the:kerneldecomp} implies the statement.
\end{proof}

The following result is an important step forward:
\begin{theorem}
  \label{the:ki}
For $\mathcal{I}=\{1,\ldots,d-1\}$ and any partition $\lambda$ with $\ell(\lambda) \le \tfrac d2$ and $s=0,\ldots,|\mathcal{I}^d_\lambda|-1$, the multiplicity of $\mathcal{H}^d_{2\lambda}$ in $H_{\lambda}^{|\lambda|+s}(\G_{d})$ is $s+1$.
\end{theorem}
\begin{proof}
For any partition $\lambda$ with $\ell(\lambda) \le \tfrac d2$, and $i=0,\dots,|\mathcal{I}^d_\lambda|-1$, there are  polynomials $q_{\lambda}^{i}:\R\to\R$ of degree $i$ satisfying the orthogonality relations
  \[
    \sum_{m=\ell(\lambda)}^{d-\ell(\lambda)} q_{\lambda}^{i}(m)q_{\lambda}^{j}(m) |v_{\lambda}^{m}|^{2} =
    \delta_{i,j}\qquad i,j=0,\dots,|\mathcal{I}^d_\lambda|-1,
  \]
  where $v_{\lambda}$ is given in \eqref{eq:Kpipi_exp}. We can define the associated positive
  definite kernels $K_{\lambda}^{i}:\G_{d}\times\G_{d}\to \R$ by
  \begin{equation}
    \label{eq:Ki}
    K_{\lambda}^{i}(P,Q) := q_{\lambda}^{i}(P)q_{\lambda}^{i}(Q) p_{\lambda}(P,Q), \qquad P,Q \in \G_{d}.
    \end{equation}
For $P \in \G_{k,d}$, $Q\in \G_{{l},d}$, the identity
    \begin{align*}
      \int_{\G_{m,d}}p_{\lambda}(P,R)p_{\lambda}(R,Q)\d\sigma_{m,d}(R) 
      &=   v_{\lambda}^{k}  |v_{\lambda}^{m}|^{2} v_{\lambda}^{{l}}
      \int_{\G_{m,d}}p_{\lambda}^{k,m}(P,R)p_{\lambda}^{m,{l}}(R,Q) \d\sigma_{m,d}(R) \\
      &= 
      |v_{\lambda}^{m}|^{2} p_{\lambda}(P,Q),
    \end{align*}
implies $K_{\lambda}^{i} * K_{\lambda}^{j} = \delta_{i,j} K_{\lambda}^{j}$. Hence, Corollary~\ref{cor:basischange} leads to
  \begin{equation}
    \label{eq:Ski_inc}
    \bigoplus_{i=0}^{s} \mathcal{S}(K^{i}_\lambda) \subset H_{\lambda}^{|\lambda|+s}(\G_{d}),\qquad
    s=0,\dots,|\mathcal{I}^d_\lambda|-1.
  \end{equation}
   Note that \eqref{eq:Ski_inc} also implies the lower bound on the multiplicities in Corollary \ref{cor:11}. 
We shall complete the proof in Appendix \ref{app:proof of long}, where we verify that \eqref{eq:Ski_inc} holds with equality, i.e.,  we decompose $H_{\lambda}^{|\lambda|+s}(\G_{d})$ into $s+1$ orthogonal subspaces of increasing polynomial degree with simple multiplicities. 
\end{proof}
We are now able to determine the multiplicities of the irreducible components in $\Pol_t(\G_{\mathcal{I},d})$.
\begin{theorem}
  \label{the:dimHkappaGKd}
  For $\mathcal{I} \subset \{1,\dots, d-1\}$ and $t\geq 0$, it holds
  \begin{equation*}
    \mu^d_{\lambda}(\mathcal{I},t) = 
    \begin{cases} 
      \min\{ t - |\lambda|+1,|\mathcal{I}^d_\lambda|\}, & 0\le|\lambda|\le t,  \quad \lambda\in\Lambda^d_\mathcal{I}, \\
      0, & \text{else}.
    \end{cases}
  \end{equation*}
\end{theorem}
\begin{proof}
According to Corollary \ref{cor:11}, we only need to suitably bound $ \mu^d_{\lambda}(\mathcal{I},t)$ from above. Proposition~\ref{the:ki} yields 
  $\mu_{\lambda}^{d}(\{1,\ldots,d-1\}^d_\lambda,t ) \leq t
  -|\lambda|+1$. Since $\mu_{\lambda}^{d}(\mathcal{I},t) \le \mu_{\lambda}^{d}(\{1,\ldots,d-1\}^d_\lambda,t )$ 
and the general upper bound
  $\mu_{\lambda}^{d}(\mathcal{I},t) \le |\mathcal{I}^d_\lambda|$ holds, we conclude the proof.
\end{proof}

The knowledge of $\dim(\mathcal{H}^d_{2\lambda})$ through \cite[Formulas
  (24.29) and (24.41)]{Fulton:1991fk} enables us to compute the dimension of $\Pol_{t}(\G_{\mathcal{I},d})$. 
Moreover, counting irreducible components yields
  \begin{equation}\label{eq:intrig 1 a}
\Pol_{t}(\G_{\mathcal{I},d}) \cong \bigoplus_{i=1}^{s}\Pol_{t-i+1}(\G_{k_{i},d}),\quad s:=\min\{t+1,|\mathcal{I}|\},
  \end{equation}
  where $\mathcal{I} =\{k_{i}\}_{i=1}^r$ is ordered such that $\min\{k_{1},d-k_{1}\} \ge \dots \ge \min\{k_{r},d-k_{r}\}$, so that there actually holds equality in Theorem~\ref{the:erstes d}. Our proof of Theorem~\ref{the:erstes d} then reveals the intriguing identity
\begin{equation}\label{eq:intrig 2}
\Null(|_{\mathcal{G}_{k_1,d}}) = (\Tr(\cdot) - k_1) \Pol_{t-1}(\G_{\mathcal{I},d}),\qquad t\geq 1,
\end{equation}
with the restriction mapping $|_{\mathcal{G}_{k_1,d}}: \Pol_{t}(\G_{\mathcal{I},d})\rightarrow  \Pol_{t}(\G_{k_1,d})$. 

\section{Cubatures and designs on unions of Grassmannians}\label{sec:cub}
So far, we have analyzed the irreducible decomposition of polynomial spaces on unions of Grassmannians. Our results enable us in the following to study \emph{cubatures} on unions of Grassmannians. 


\subsection{Introducing cubatures and designs}\label{sec:cub design}
Any orthogonally invariant finite signed measure $\sigma_{\mathcal{I},d}$ on
$\mathcal{G}_{\mathcal{I},d}$ is a linear combination of the Haar
(probability) measures $\sigma_{k,d}$, $k \in\mathcal{I}$, i.e.,
\begin{equation*}
  \sigma_{\mathcal{I},d} = \sum_{k\in\mathcal{I}}  m_k \sigma_{k,d},
\end{equation*}
for some $\{m_k\}_{k\in\mathcal{I}}\subset\R$. For
points $\{P_j\}_{j=1}^n\subset \mathcal{G}_{\mathcal{I},d}$ and
 weights $\{\omega_j\}_{j=1}^n \subset \R$, we say that $\{(P_j,\omega_j)\}_{j=1}^n$ is a cubature for
$\Pol_{t}(\mathcal{G}_{\mathcal{I},d})$ (resp.~$\Hom_{t}(\G_{\mathcal{I},d})$) with
respect to $\sigma_{\mathcal{I},d}$ if
\begin{equation*}
  \int_{\mathcal{G}_{\mathcal{I},d}} f(P)d\sigma_{\mathcal{I},d}(P) = \sum_{j=1}^n \omega_j f(P_j),\qquad f\in  \Pol_{t}(\mathcal{G}_{\mathcal{I},d})\quad (\mathrm{resp.}\; f\in\Hom_{t}(\G_{\mathcal{I},d})).
\end{equation*}
If $\{(P_j,\omega_j)\}_{j=1}^n$ is a cubature for
$\Hom_{t}(\mathcal{G}_{\mathcal{I},d})$ with
$\sum_{k \in \mathcal{I}} m_{k} = \sum_{j=1}^{n} \omega_{j}$, then it is also a cubature
for $\Pol_{t}(\G_{\mathcal{I},d})$, cf.~Theorem~\ref{the:structure}. The value of the parameter $t$ is often called the \emph{strength} of the cubature. 
\begin{remark}\label{rem:pache}
According to \cite[Proposition~2.6 and 2.7]{Harpe:2005fk}, 
there is a cubature $\{(P_j,\omega_j)\}_{j=1}^n$ with nonnegative weights $\{\omega_j\}_{j=1}^n$ for $\Pol_{ t}(\mathcal{G}_{\mathcal{I},d})$ with respect to $\sigma_{\mathcal{I},d}$ such that $n\leq  \dim(\Pol_{ t}(\mathcal{G}_{\mathcal{I},d}))$. For $m_k\neq 0$, $k\in\mathcal I$, any cubature $\{(P_j,\omega_j)\}_{j=1}^n$ for $\Pol_{2t}(\mathcal{G}_{\mathcal{I},d})$ with respect to $\sigma_{\mathcal{I},d}$ with nonnegative weights satisfies $n\geq \dim( \Pol_{t}(\mathcal{G}_{\mathcal{I},d}))$, cf.~\cite[Proposition~1.7]{Harpe:2005fk}.
\end{remark}

Analogous to Euclidean designs, cf.~\cite{Bajnok:2006jt,Bajnok:2007wq,Bannai:2006la,Bannai:2012wd}, cubatures on unions of Grassmannians induce cubatures on single Grassmannians, but potentially with lower strength: 
\begin{proposition}\label{prop:marginal designs}
  If $\{(P_j,\omega_j)\}_{j=1}^n$ is a cubature for
  $\Pol_{t}(\mathcal{G}_{\mathcal{I},d})$, $t\ge |\mathcal{I}|-1$, with respect
  to the signed measure
  $\sigma_{\mathcal{I},d}=\sum_{k\in\mathcal{I}} m_k \sigma_{k,d}$, then, for
  any $k\in\mathcal{I}$,
\begin{equation*}
  \{ (P_j,\omega_j) : P_j\in\mathcal{G}_{k,d},\; j=1,\ldots,n\}
\end{equation*}
is a cubature for $\Pol_{s}(\mathcal{G}_{k,d})$ with respect to the signed
measure $m_{k} \sigma_{k,d}$, where $s=t-|\mathcal{I}|+1$.
\end{proposition}
\begin{proof}
  Let $f$ be a polynomial of degree at most $s$ on $\R^{d\times d}_{\sym}$, then we know
  that $1_{\mathcal{G}_{k,d}} \!\cdot f$ is a polynomial of degree at most $t$.
  Hence, the statement follows from
\[
  m_{k} \int_{\mathcal{G}_{k,d}} f(P)d\sigma_{k,d} = \int_{\mathcal{G}_{\mathcal{I},d}} (1_{\mathcal{G}_{k,d}}\! \cdot f)(P)d\sigma_{\mathcal{I},d} 
  = \sum_{\{j:P_j\in\mathcal{G}_{k,d}\}} \omega_j f(P_j).\qedhere
\]
\end{proof}

We also observe that any cubature of strength $2t$ gives rise to a cubature of strength $2t+1$:
\begin{proposition}\label{prop:2t to 2t+1}
Let $\{(P_j,\omega_j)\}_{j=1}^n$ be a cubature for $\Pol_{t}(\mathcal{G}_{\mathcal{I},d})$ with respect to the signed measure $\sigma_{\mathcal{I},d}=\sum_{k\in\mathcal{I}} m_k \sigma_{k,d}$.  Then 
\begin{equation*}
\{(P_j,\omega_j)\}_{j=1}^n \cup \{(I-P_j,(-1)^t\omega_j)\}_{j=1}^n
\end{equation*}
is a cubature for $\Pol_{t+1}(\mathcal{G}_{\mathcal{I}\cup (d-\mathcal{I}),d})$ with respect to the signed measure  
\begin{equation*}
  \sigma_{\mathcal{I}\cup(d-\mathcal{I}),d}=\sum_{k\in\mathcal{I}} m_k(\sigma_{k,d} +(-1)^t\sigma_{d-k,d}).
\end{equation*}
\end{proposition}
\begin{proof}
As in the proof of Theorem~\ref{the:structure}, it is sufficient to consider
the polynomial  $\Tr( X \cdot)^{t+1}|_{\G_{\mathcal{I}\cup(d-\mathcal{I}),d}}$ 
  for $X\in\R^{d\times d}_{\sym}$. Then we have
\begin{align*}
    \int_{\mathcal{G}_{\mathcal{I}\cup (d-\mathcal{I}),d} }\!\!\!\! \!\!\!\!\!\!\!\!\!\!\!\!\!\!\!\! \Tr( XP)^{t+1} d\sigma_{\mathcal{I}\cup(d-\mathcal{I}),d}(P) &= 
    \int_{\mathcal{G}_{\mathcal{I},d}} \!\!\!\!\!\!\Tr( XP)^{t+1}d\sigma_{\mathcal{I},d} (P) + \!\!
 \int_{\mathcal{G}_{\mathcal{I},d}} \!\!\!\!\!\!\Tr( X(I-P))^{t+1}(-1)^td\sigma_{\mathcal{I},d} (P) \\
 &    =\int_{\mathcal{G}_{\mathcal{I},d}} \!\!\!\!\!\!\big(\Tr(XP)^{t+1} +
    (-1)^t\Tr(X(I-P))^{t+1} \big)d\sigma_{\mathcal{I},d}(P).
\end{align*}
The mapping $P\mapsto \Tr( XP)^{t+1} +  (-1)^t\Tr( X(I-P))^{t+1} $ is
contained in $\Pol_{t}(\mathcal{G}_{\mathcal{I},d})$, because the two terms with exponent $t+1$ cancel out. Thus,
$\{(P_j,\omega_j)\}_{j=1}^n$ being a cubature for
$\Pol_{t}(\mathcal{G}_{\mathcal{I},d})$ yields
\[
\int_{\mathcal{G}_{\mathcal{I}\cup (d-\mathcal{I}),d} }\!\!\!\!\!\! \Tr(
XP)^{t+1} d\sigma_{\mathcal{I}\cup(d-\mathcal{I}),d}(P) = \sum_{j=1}^n
\omega_j \big( \Tr(XP_j)^{t+1} + (-1)^t\Tr( X(I-P_j))^{t+1}
\big).\qedhere
\]
\end{proof}


Cubatures of strength $t$ on single Grassmannians, whose weights are all the same, are called $t$-designs and have been studied in \cite{Bachoc:2005aa,Bachoc:2006aa,Bachoc:2004fk,Bachoc:2002aa}. We shall extend this concept to unions of Grassmannians after a brief observation on characterstic functions. For $t=|\mathcal{I}|-1$ with $\mathcal{I}=\{k_0,\ldots,k_{t}\}$, the Vandermonde matrix $V=(k_i^j)_{i,j=0,\ldots,t}\in\R^{t+1}$ is invertible. Hence, for $i=0,\ldots,t$, there are $\alpha_i:=(\alpha_{i,0},\dots,\alpha_{i,t})^\top \in\R^{t+1}$ with $V
\alpha_i=e_{i+1}$ implying 
\begin{equation*}
\sum_{j=0}^{t} \alpha_{i,j} \Tr(P)^j=\begin{cases}
1,&  P\in\mathcal{G}_{k_i,d},\\
0,& P\in\mathcal{G}_{\mathcal{I},d}\setminus \mathcal{G}_{k_i,d}.
\end{cases}
\end{equation*}
Thus, the characteristic function of each $\mathcal{G}_{k,d}$ is contained in
$\Pol_t(\mathcal{G}_{\mathcal{I},d})$, for all $t\geq |\mathcal{I}|-1$. The latter yields that any cubature $\{(P_j,\omega_j)\}_{j=1}^n$ for $\Pol_t(\mathcal{G}_{\mathcal{I},d})$ with respect to 
$\sigma_{\mathcal{I},d}=\sum_{k\in\mathcal{I}}m_k \sigma_{k,d}$ satisfies
$ m_k=\sum_{\{j:P_j\in\mathcal{G}_{k,d}\}} \omega_j$, for $k\in \mathcal{I}$. 
If weights are the same on each single $\mathcal{G}_{k,d}$, then those weights must be $\frac{m_k}{n_k}$,
where $n_k=|\{j:P_j\in\mathcal{G}_{k,d}\}|$. Indeed, we impose 
this condition in our definition of designs:
\begin{definition}
Let $\mathcal{P}_k\subset\mathcal{G}_{k,d}$ be finite and denote $n_k:=|\mathcal{P}_k|$. 
The collection $\mathcal{P}_\mathcal{I}=\bigcup_{k\in\mathcal{I}}\mathcal{P}_k$ is called a $t$-design with respect to an orthogonally invariant signed measure $\sigma_{\mathcal{I},d}=\sum_{k\in\mathcal{I}}m_k \sigma_{k,d}$ if
\begin{equation*}
\int_{\mathcal{G}_{\mathcal{I},d}} f(P)d\sigma_{\mathcal{I},d}(P) = \sum_{k\in\mathcal{I}} \tfrac{m_k}{n_k} \sum_{P\in\mathcal{P}_k} f(P),\quad\text{for all } f\in\Pol_t(\mathcal{G}_{\mathcal{I},d}).
\end{equation*}
\end{definition}
Thus, any $t$-design is a cubature, whose weights are the same on each Grassmannian but can differ across different Grassmannians. 
For $|\mathcal{I}|=1$, our definition reduces to the standard Grassmannian designs as considered in \cite{Bachoc:2005aa,Bachoc:2006aa,Bachoc:2004fk,Bachoc:2002aa}. The existence of $t$-designs with $t=1$ in single Grassmannians was  studied in \cite{Casazza:2011aa}. For a discussion on the existence of
  $1$-designs in unions of Grassmannians with
  $\frac{m_k}{n_k}=\frac{m_{l}}{n_{l}}$, for all $k,{l}\in\mathcal{I}$, we refer to
  \cite{Bownik:la}. 

\subsection{Constructing cubatures and designs by numerical minimization}\label{sec:potential}
To construct cubatures or designs, we consider the 
$t$-fusion frame potential, cf.~\cite{Bachoc:2010aa},
\begin{equation*}
  \FFP_t(\{(P_j,\omega_j)\}_{j=1}^n):=\sum_{i,j=1}^{n}\omega_i\omega_j\Tr(P_iP_j)^t.
\end{equation*}
where $\{P_j\}_{j=1}^n\subset\mathcal{G}_{\mathcal{I},d}$ and $\{\omega_j\}_{j=1}^n$. The 
$1$-fusion frame potential was already investigated in \cite{Casazza:2009aa,Massey:2010fk}. Lower bounds on the
$t$-fusion frame potential for general positive integers $t$ were derived in
\cite{Bachoc:2010aa}. Only for single Grassmannians, i.e., $|\mathcal{I}|=1$, those
$\{(P_j,\omega_j)\}_{j=1}^n$ were characterized in \cite{Bachoc:2010aa}, for which the bounds are matched. In the following, we provide a characterization for the general case $|\mathcal{I}| \geq 1$. Before we state this result though, it is convenient to define 
\begin{equation}
  \label{eq:Tsigma}
  \mathcal{T}_{\sigma_{\mathcal{I},d}}(t):=\int_{\mathcal{G}_{\mathcal{I},d}}\int_{\mathcal{G}_{\mathcal{I},d}} \Tr( PQ)^td\sigma_{\mathcal{I},d}(P)d\sigma_{\mathcal{I},d}(Q) = m^{\top} T_{\mathcal{I},d}(t) m,
\end{equation}
where $\sigma_{\mathcal{I},d}=\sum_{k\in\mathcal{I}}m_k \sigma_{k,d}$ and $m=(m_{k})_{k\in \mathcal{I}}$ with the matrix
$T_{\mathcal{I},d}(t) \in \R^{|\mathcal{I}|\times |\mathcal{I}|}$ being given by 
\begin{equation*}
(T_{\mathcal{I},d}(t))_{k,{l}}:= \int_{\G_{{l},d}}\int_{\G_{k,d}}\Tr(P Q)^{t}
d\sigma_{k,d}(P) d\sigma_{{l},d}(Q). 
\end{equation*}
Note that $T_{\mathcal{I},d}(t)$ is the $(0)$-th Fourier coefficient of the
positive definite zonal kernel $K_{t}(P,Q) = \Tr( P Q)^{t}$, for $P,Q \in \G_{\mathcal{I},d}$, and thus symmetric and positive
semidefinite.

\begin{theorem}\label{the:pot}
  Given $\{P_j\}_{j=1}^n\subset \mathcal{G}_{\mathcal{I},d}$ with weights $\{\omega_j\}_{j=1}^n \subset \R$, let $m_{k}=\sum_{\{j:P_j\in\mathcal{G}_{k,d}\}}\omega_j$. 
  Then the fusion frame potential is bounded from below by
 \begin{equation}\label{eq:ene func}
    \FFP_t(\{(P_j,\omega_j)\}_{j=1}^n)\geq  \mathcal{T}_{\sigma_{\mathcal{I},d}}(t).
  \end{equation}
  Equality holds if and only if
  $\{(P_j,\omega_j)\}_{j=1}^n$ is a cubature for $\Pol_{
    t}(\mathcal{G}_{\mathcal{I},d})$ with respect to $\sigma_{\mathcal{I},d}$.
%
\end{theorem}
This theorem is an extension of results in \cite{Bachoc:2010aa}, where the lower bound is already
derived but equality is not discussed. For the proof, we refer to Appendix \ref{sec:proof 6.5}.


Theorem~\ref{the:pot} enables the use of numerical minimization schemes to derive cubatures. Knowledge of the global lower bound $\mathcal{T}_{\sigma_{\mathcal{I},d}}(t)$ is important to check if numerical solutions are indeed cubatures by ruling out that the minimization got stuck in a local minimum. The matrix
$T_{\mathcal{I},d}(t)$ can be computed via zonal
polynomials
by 
\begin{equation*}
  T_{\mathcal{I},d}(t) = \sum_{\substack{
      |\lambda| = t,\\
      \ell(\lambda) \le d/2}
  } \frac{c_{\lambda}
    c_{\lambda}^{\top}}{C_{\lambda}(I_{d})} \in \mathbb R^{|\mathcal{I}| \times |\mathcal{I}|},
  \qquad c_{\lambda}:=(C_{\lambda}(I_{k}))_{k \in \mathcal{I}} \in \mathbb R^{|\mathcal{I}|},
\end{equation*}
cf.~\eqref{eq:intCkappa} in Appendix \ref{sec:zon}, and note that $C_{\lambda}(I_k)$ is explicitly computed in \cite{Chikuse:2003aa,Muirhead:1982fk,Gross:1987bf}.

For suitable minimization algorithms on Grassmannians, we refer to \cite{Graf:2013zl} and
  \cite{Absil:2008qr,Boumal:2014ff}. In Section \ref{sec:examples num} we shall indeed numerically minimize the fusion frame potential with equal weights and check that the lower bound is attained. 


\section{Examples of $t$-designs derived from numerical minimization}\label{sec:examples num}
Here, we shall construct some families of $1$-, $2$-, and $3$-designs in unions of Grassmannians. By numerically minimizing the energy functional \eqref{eq:ene func} using a conjugate gradient approach, cf.~\cite[Section 3.3.1]{Graf:2013zl}, we compute candidates for $t$-designs, i.e., $t$-designs up to machine precision. Based on the special structures of the Gram matrices of these numerical minimizers, we looked for group orbits that describe them analytically. Indeed ``beautifying'' our numerical results, we were able to analytically specify our candidates, which turned out to be exact minimizers.

\subsection{A family of $1$-designs in arbitrary dimensions}\label{sec:9.11}
We analytically construct $d$ lines and $1$ hyperplane in $\R^d$, so that the corresponding orthogonal projectors are a $1$-design in $\mathcal{G}_{1,d}\cup\mathcal{G}_{d-1,d}$ with respect to the measure $\sigma_{\{1,d-1\},d}=m_{1}\sigma_{1,d}+m_{d-1}\sigma_{d-1,d}$, with $m_{1}=1$ and $m_{d-1}\in [-1/(d-1),1]$. Indeed, let the $d$ lines be described by the vectors
\[
a_{i} := \sqrt{1-m_{d-1}}\, e_{i} + \frac{1}{d}\big(\sqrt{1+(d-1)m_{d-1}}-\sqrt{1-m_{d-1}} \big) e\, \in \mathbb R^{d}, \qquad i=1,\dots,d,
\]
where $e_{i} \in \mathbb R^{d}$ are the standard unit vectors, and
$e:=e_{1}+\dots+e_{d}$. The associated rank-1 projectors are $P_{i}:=a_{i}\,
a_{i}^{\top} \in \mathbb R^{d}$, $i=1,\dots,d,$ and the rank-$(d-1)$ projector
is $P_{d+1}:=I_{d} - \frac{1}{d} e\, e^{\top}$. We calculate
\[
\mathcal T_{\sigma_{1,d}+m_{d-1}\sigma_{{d-1,d}}}(1) 
= \tfrac{1}{d}\big( 1 + (d-1) m_{d-1}\big)^{2}
\]
and observe equality with the $1$-fusion frame potential
$\FFP_1\big(\{(P_{j},\omega_{j})\}_{j=1}^{d+1}\big)$ with weights
$\omega_{1}=\dots=\omega_{d}=\frac{m_{1}}{d}$, $\omega_{d+1}=m_{d-1}$. According to Theorem~\ref{the:pot}, $\{(P_{j},\omega_{j})\}_{j=1}^{d+1}$ is a $1$-design. 

For $m_{d-1}=1$, the $d$ lines coincide and are orthogonal to the hyperplane. This configuration can also be considered as a single line with one hyperplane forming a $1$-design. If $-1/(d-1)\leq m_{d-1}<1$, then this one line is split into $d$ lines forming a $(d-1)$-dimensional simplex, and, for $m_{d-1}=-1/(d-1)$, the lines lie in the hyperplane. 


\begin{remark}
  According to Proposition~\ref{prop:marginal designs}, any cubature of strength $1$
  in $\mathcal{G}_{1,d}\cup \mathcal{G}_{d-1,d}$ yields parts in
  $\mathcal{G}_{1,d}$ and in $\mathcal{G}_{d-1,d}$ that are cubatures of
  strength $0$. Hence, there must at least be one element in each. For $m_{1}=m_{d-1}=1$, the above example matches this lower bound. 
\end{remark}

\subsection{A family of $1$-designs in $\R^4$} 
We construct a family of $2$ lines and $2$ planes in $\R^4$ forming a $1$-design
with respect to $m_1\in [-2,2]$ and $m_2=1$. The rank-1 projectors are  $
P_{1}:=e_{1}\,e_{1}^{\top}$, $P_{2}:=e_{3}\,e_{3}^{\top} \in \mathbb R^{4}$, 
and the $2$-dimensional projectors are
\begin{equation*}
  P_{3} := 
  \begin{pmatrix}
    a_{3}\,a_{3}^{\top} & 0 \\
    0 & a_{3}\,a_{3}^{\top}
  \end{pmatrix}, \;
  P_{4} := 
  \begin{pmatrix}
    a_{4}\,a_{4}^{\top} & 0 \\
    0 & a_{4}\,a_{4}^{\top}
  \end{pmatrix} \in \mathbb R^{4}, \quad a_{3/4}:=\frac12
  \begin{pmatrix}
    \sqrt{2 - m_{1}}\\
    \pm \sqrt{2 + m_{1}}
  \end{pmatrix}
  \in \R^{2}.
\end{equation*}
For $m_{1} \in [-2,2]$, this family provides $1$-designs since there holds equality in Theorem \ref{the:pot}, where 
$ 
\mathcal T_{m_{1} \sigma_{1,4}+ \sigma_{2,4}}(1)=\left(1 + \tfrac{m_{1}}{2}
\right)^{2}.
$ 
For $m_1=2$, the two planes coincide and the two lines are orthogonal to each
other and to the planes. The two planes also
coincide for $m_1=-2$ and then the two lines span the same plane. The choice
$m_1=0$ yields two planes that are orthogonal to each other.
\begin{remark}
For $m_1,m_2>0$, there do not exist any $1$-designs in $\mathcal{G}_{1,4}\cup \mathcal{G}_{2,4}$ of cardinality $2$. In this sense, our example with $m_1=2$ and $m_2=1$
  is optimal.
\end{remark}

\subsection{A family of $2$-designs in $\R^3$} \label{sec:2 designs in
  R3} 
The numerical minimization enabled us to identify an analytic family of $6$ lines and $4$ planes in $\R^3$ forming a $2$-design with $m_1=1$ and $m_2\in[-3/8,3/2]$. The tetrahedral group $\mathrm T$ is 
generated by a cyclic coordinate shift $R_{1}=e_2e_1^\top+e_3e_2^\top+e_1e_3^\top\in\R^{3\times 3}$ and a reflection $R_{2}=I_3-2e_1e_1^\top$. 
The 6 lines and 4 planes are obtained as distinct orbits under the canonical action of the
tetrahedral group. The rank-1 and rank-2 projectors are
\[
\begin{aligned}
  \mathcal P_{1} & := \{ R P_{1} R^{\top} \;:\; R \in \mathrm T \}, \qquad
  P_{1}:=\frac12 a_{1}\,a_{1}^{\top}, \quad a_{1}:=\Big(\sqrt{1-v},\sqrt{1+v},0\Big)^{\top},\\
  \mathcal P_{2} & := \{ R P_{2} R^{\top} \;:\; R \in \mathrm T \}, \qquad
  P_{2}:=I - \frac13 e\,e^{\top} \in \mathbb R^{3},
\end{aligned}
\]
respectively, where $v:=\sqrt{\frac{1}{15}(3+8m_{2})}$. For $m_{2}\in [-3/8,3/2]$, the union 
$\mathcal P=\mathcal P_{1} \cup \mathcal P_{2}$ is a 2-design since there holds equality in Theorem \ref{the:pot} for $t=2$ with 
$ 
\mathcal T_{\sigma_{1,3}+ m_{2}
  \sigma_{2,3}}(2)=\tfrac{1}{15}\left(3+16m_{2}+28m_{2}^{2}\right).
$ 
For $m_2=3/2$, there are $4$ planes and the $6$ lines collapse to the $3$ coordinate axes. 
If $m_2=-3/8$, then the $6$ lines are exactly the $6$ intersection lines of the $4$ planes. For $m_2=0$, the $6$ lines correspond to the vertices of the icosahedron. 

\begin{remark}
  Any cubature of strength
  $2$ with nonnegative weights needs at least $n$ distinct cubature points with
  $n\geq \dim(\Pol_1(\mathcal{G}_{\mathcal{I},d}))$, cf.~Remark~\ref{rem:pache}. Theorem~\ref{the:dimHkappaGKd} leads to
  $\dim(\Pol_1(\mathcal{G}_{1,3} \cup \mathcal{G}_{2,3})) = 7$, so that our construction of $4$ planes and $3$ lines
  associated to the choice $m_1=1$ and $m_2=3/2$ has the minimal number of
  elements. Moreover, Proposition~\ref{prop:marginal designs} implies that any $2$-design on
  $\mathcal{G}_{1,3} \cup \mathcal{G}_{2,3}$ has at least $3$ lines.
\end{remark}

\subsection{A family of $2$-designs in $\R^4$} \label{sec:2 designs in
  R4} We shall provide a family of $8$ lines and $8$ planes in $\R^4$ forming a
$2$-design for $m_1=1$ and $m_2\in[3/4,3/2]$. 
The lines are generated by a symmetry group $\mathrm G_{1}$ of cardinality $|\mathrm G_{1}|= 64$, which is generated by a cyclic coordinate shift $R_{1}=e_2e_1^\top+e_3e_2^\top+e_4e_3^\top+e_1e_4^\top\in\R^{4\times 4}$ and
reflection $R_{2}=I_4-2e_1e_1^\top$. 
The rank-1 projectors are the orbit
\[
\mathcal P_{1}:= \{ R P_{1} R^{\top} \;:\; R \in \mathrm G_{1} \}, \quad
P_{1}:=\frac12 a_{1}\,a_{1}^{\top}, \;\;
a_{1}:=\Big(\sqrt{1-v},\sqrt{1+v},0,0\Big)^{\top},
\]
where $v:=\sqrt{\frac19(8m_{2}-3)}$. The planes are derived from two orbits
under a cyclic group $\mathrm G_{2}$ of cardinality $|\mathrm G_{2}|=4$
generated by $R_3=-e_3e_1^\top+e_2e_2^\top+e_1e_3^\top-e_4e_4^\top\in\R^{4\times 4}$. 
More precisely, the set of rank-2 projectors is $\mathcal P_{2}  := \{ R P_{i} R^{\top} \;:\; R \in \mathrm G_{2},\; i=2,3 \}$ with 
\[
\begin{aligned}
P_{2/3} & := \frac16
\begin{pmatrix}
 3\mp\sqrt{1-b^{2}}             &     \sqrt{2-3b+b^{2}} &   \pm\sqrt{4-b^{2}}             & \mp\sqrt{2+3b+b^{2}} \\ 
  \sqrt{2-3b+b^{2}}    & 3\pm\sqrt{1-b^{2}} &  \sqrt{2+3b+b^{2}}    & \sqrt{4-b^{2}} \\ 
  \pm\sqrt{4-b^{2}}                &     \sqrt{2+3b+b^{2}} &  3\mp\sqrt{1-b^{2}}           & \mp\sqrt{2-3b+b^{2}} \\
  \mp\sqrt{2+3b+b^{2}} &       \sqrt{4-b^{2}}      &  \mp\sqrt{2-3b+b^{2}}     & 3\pm\sqrt{1-b^{2}} 
\end{pmatrix},
\end{aligned}
\]
where $b:=\sqrt{2-\frac{3}{2m_{2}}}$. For $\mathcal P:=\mathcal P_{1}\cup
\mathcal P_{2}$, the lower bound 
$ 
\mathcal
T_{\sigma_{1,4}+m_{2}\sigma_{2,4}}(2)=\tfrac{1}{72} \left(9 + 48 m_{2} + 80 m_{2}^2\right)
$ 
 on the $2$-fusion frame potential is matched for $m_{2}\in [3/4,3/2]$. Therefore, $\mathcal P$ is
indeed a 2-design. For $m_2=3/2$, each plane intersects two other planes in one line, and the $8$
lines coincide in pairs with the $4$ coordinate axes.

\begin{remark}
Proposition~\ref{prop:marginal designs} yields that any $2$-designs with $d=4$, $m_1=1$ and $m_2=3/2$ must have at least $4$ lines, which is matched by our example of $4$ lines and $8$ planes.
\end{remark}

\subsection{A $2$-design in $\R^5$} \label{sec:2 design in R5} The following
collection of $5$ lines and $16$ planes in $\R^5$ forms a $2$-design for $m_1=1$
and $m_2=5/3$. The lines are given by the $5$ coordinate axes, i.e., 
$\mathcal P_{1}:=\{ P_{i}: i=1,\dots,5 \}$, $P_{i}:= e_{i} e_{i}^{\top}
\in \R^{5\times 5}$, for $i=1,\dots,5$. To construct the planes, we consider the transformation group
$\mathrm G$ of cardinality $|\mathrm G| = 16$, generated by $G_1=I_3 - e_5e_4^\top+e_4e_5^\top \in\R^{5\times 5}$ and $G_2=e_2e_1^\top-e_1e_2^\top+e_3e_3^\top+e_4e_4^\top-e_5e_5^\top$.  
The set of rank-2 projectors is the orbit $\mathcal P_{2} :=\{ R P_{6} R^{\top} \;:\; R \in \mathrm G \}$, where 
\begin{equation*}
P_{6}  := \frac15 \left(
  \begin{smallmatrix}
2 &\sqrt {\frac {3} {2}} &\sqrt {\frac {1} {6}\left (9 + 
                    5\sqrt {3} \right)} &\frac {1} {6}\left (3 - 
                    5\sqrt {3} \right) &\sqrt {\frac {2} {3}} 
\\\sqrt {\frac {3} {2}} & 2 &\sqrt {\frac {1} {6}\left (9 - 
                    5\sqrt {3} \right)} &\sqrt {\frac {2} {3}} 
&\frac {1} {6}\left (3 + 
                    5\sqrt {3} \right) \\\sqrt {\frac {1} {6}\left 
(9 + 5\sqrt {3} \right)} &\sqrt {\frac {1} {6}\left (9 - 
                    5\sqrt {3} \right)} & 2 & - \sqrt {\frac {1} {6}
\left (9 + 5\sqrt {3} \right)} & - \sqrt {\frac {1} {6}\left (9 - 
                    5\sqrt {3} \right)} \\\frac {1} {6}\left (3 - 
                    5\sqrt {3} \right) &\sqrt {\frac {2} {3}} & - 
\sqrt {\frac {1} {6}\left (9 + 
                    5\sqrt {3} \right)} & 2 &\sqrt {\frac {3} {2}} 
\\\sqrt {\frac {2} {3}} &\frac {1} {6}\left (3 + 
          5\sqrt {3} \right) & - \sqrt {\frac {1} {6}\left (9 - 
           5\sqrt {3} \right)} &\sqrt {\frac {3} {2}} & 2 
  \end{smallmatrix} 
\right).
\end{equation*}
For $\mathcal
P:=\mathcal P_{1}\cup \mathcal P_{2}$, the lower bound 
$ 
\mathcal T_{\sigma_{1,5}+\frac53
  \sigma_{2,5}}(2)=\tfrac{131}{45}
$ 
on the $2$-fusion frame potential is matched, so that $\mathcal P$ is indeed a
2-design.

\begin{remark}
Similar to the previous example, Proposition~\ref{prop:marginal designs} yields that any $2$-design with $d=5$, $m_1=1$ and $m_2=5/3$ must have at least $5$ lines, which is matched.
\end{remark}

\subsection{From $2t$-designs to $2t+1$-designs}
Theorem~\ref{prop:2t to 2t+1} yields a construction of $2t+1$-designs from
$2t$-designs. Hence, any $2$-design in the Sections \ref{sec:2 designs in R3}, \ref{sec:2 designs in R4}, and \ref{sec:2 design in R5} gives rise to a $3$-design. For instance, the $6$ lines going through the vertices of the icosahedron are a $2$-design in $\R^3$. By adding the $6$ complementary planes, we obtain a $3$-design with respect to $m_1=m_2=1$. 
\begin{remark}
Proposition~\ref{prop:marginal designs} yields that the parts in each single Grassmannian $\mathcal{G}_{k,d}$ of a $3$-design for $\mathcal{G}_{\mathcal{I},d}$ with $|\mathcal{I}|=2$ must be a cubature of strength $2$. Any such cubature must have at least $6$ elements for $d=3$, so that our $3$-design with $12$ elements in the union of two Grassmannians is optimal. 
\end{remark}

%

\section*{Acknowledgements}
The authors would like to thank Christine Bachoc and Christian Krattenthaler for valuable advice and discussions. Both authors have been funded by the Vienna Science and Technology Fund (WWTF) through project VRG12-009.

\bibliographystyle{amsplain}
\bibliography{../biblio_ehler2}

\newpage

\begin{appendix}

\section{Reproducing kernels}
\subsection{Some basics on reproducing kernels}\label{app:rep k}
A Hilbert space $\mathscr{H}(X)$ of continuous functions on a set $X\subset\R^d$ with inner product
$(f,g)_{\mathscr{H}(X)}$, for $f,g \in \mathscr{H}(X)$, is called a \emph{reproducing kernel Hilbert space} if
point evaluation is continuous. For such a Hilbert space there exists a unique (continuous) positive definite kernel $K:X\times X\to \C$, i.e., $K(x,y) = \overline{K(y,x)}$, for $x,y \in X $, and 
\begin{equation*}
0\leq \sum_{i,j=1}^{M} c_{i} \overline{c_{j}} K(x_{i},x_{j}), \quad\text{for all } x_{i}
  \in X,\; c_{i} \in \C,\; i=1,\dots,M,\; M \in \N,
\end{equation*}
such that $K(x,\cdot) \in \mathscr{H}(X)$, $x \in X$, and the reproducing
property
\begin{equation}\label{eq:rep def}
f(x) = (f, K(x,\cdot))_{\mathscr{H}(X)},\quad f \in \mathscr{H}(X), \quad x \in X,
\end{equation}
holds. 

Conversely, any continuous positive definite kernel
$K:X\times X\to \C$ gives rise to a unique reproducing kernel Hilbert
space $\mathcal{S}(K)$ of continuous functions, so that the function space
\begin{equation}\label{eq:S0 def}
\mathcal{S}^0(K) := \spann\{ K(x,\cdot) :x \in X \}
\end{equation}
is dense in $\mathcal{S}(K)$ and $K$ satisfies the reproducing property. Here, we exclusively deal with polynomial kernels $K$ so that $m:=\dim(\mathcal{S}^0(K))<\infty$, and a basis of $\mathcal{S}(K)=\mathcal{S}^0(K)$ is given by the functions $
K(x_{i},\cdot):X\to \C$, for $i=1,\dots,m$, 
if the matrix
\begin{equation*}
\big(K(x_{i},x_{j})\big)_{i,j=1}^{m} \in \C^{m\times m}
\end{equation*}
is invertible. The reproducing property \eqref{eq:rep def} is satisfied with respect to the inner product 
\begin{equation*}
(f,g)_{\mathcal{S}(K)} := \sum_{i,j=1}^m \alpha_i\beta_jK(x_i,x_j),
\end{equation*}
where $f=\sum_{i=1}^m\alpha_i K(x_i,\cdot)$ and $g=\sum_{i=1}^m\beta_i K(x_i,\cdot)$. It turns out that the restrictions $\mathcal{S}(K)|_Y$ to compact subsets $Y\subset X$ are generated by the restricted kernel itself, i.e., if $Y \subset X$ is compact, then
\begin{equation} \label{eq:HKYspann}
\mathcal{S}(K)|_Y = \mathcal{S}(K|_{Y\times Y}),
\end{equation}
see also \cite[Theorem~10.47]{Wendland:2004wd} for a more general setting. 
For more details on the theory of reproducing kernel Hilbert spaces, we refer to  \cite{Wendland:2004wd}.

%

\subsection{Proof of Theorem~\ref{the:pot}}\label{sec:proof 6.5}
For the sake of completeness, we provide the proof of Theorem~\ref{the:pot}. 
\begin{proof}
The worst-case cubature error $E_{t}$ in the reproducing kernel Hilbert space
$\mathcal{S}(K_t) = \Hom_{t}(\G_{\mathcal{I},d})$ with reproducing kernel
$K_{t}(P,Q) = \Tr( P Q)^{t}$, for $P,Q \in \G_{\mathcal{I},d}$,  is defined by
\begin{equation}
  \label{eq:wcce}
  E_{t}(\{(P_{j},\omega_{j})\}_{j=1}^{n},\sigma_{\mathcal{I},d}) := \sup_{f\in \Hom_{
      t}(\mathcal{G}_{\mathcal{I},d}), \| f \|_{K_{t}} = 1} \Big| \int_{\G_{\mathcal{I},d}}
  f(P) \d\sigma_{\mathcal{I},d}(P) - \sum_{i=1}^{n} \omega_{i} f(P_{i})\Big|,
\end{equation}
where $\| \cdot\|_{K_{t}}$ is the norm in $\Hom_{t}(\G_{\mathcal{I},d})$ induced by
$K_{t}$. The squared worst-case cubature error satisfies
\begin{align*}
    E_{t}^{2}(\{(P_{j},\omega_{j})\}_{j=1}^{n},\sigma_{\mathcal{I},d}) 
    & = \sum_{i,j=1}^{n}\omega_{i}\omega_{j}
    \Tr(P_{i}P_{j})^t - 2 \sum_{i=1}^{n} \omega_{i} \int_{\G_{\mathcal
        K,d}} \Tr(P_{i}Q)^t \d\sigma_{\mathcal{I},d}(Q) \\
    & \qquad + \int_{\G_{\mathcal
        K,d}}\int_{\G_{\mathcal{I},d}} \Tr(PQ)^t \d\sigma_{\mathcal{I},d}(P)\d\sigma_{\mathcal{I},d}(Q),
\end{align*}
cf.~\cite[Theorem~2.7]{Graf:2013zl}. Isolating the $t$-fusion frame potential and noting that the conditions on $m_k$ imply $\sum_{i=1}^{n} \omega_{i} \int_{\G_{\mathcal
        K,d}} K_{t}(P_{i},Q) \d\sigma_{\mathcal{I},d}(Q)=  \mathcal{T}_{\sigma_{\mathcal{I},d}}(t)$ yields the lower bound. The expression in \eqref{eq:wcce} equals $0$ if
and only if $\{(P_{j},\omega_{j})\}_{j=1}^{n}$ is a cubature for
$\Hom_{ t}(\mathcal{G}_{\mathcal{I},d})$ with respect to $\sigma_{\mathcal{I},d}$. 
\end{proof}

\section{Properties of special polynomials}
\subsection{Zonal polynomials}\label{sec:zon}
The zonal polynomial $C_\lambda$ is an orthogonally invariant homogeneous polynomial of degree $|\lambda|$ on
$\R^{d\times d}_{\sym}$, see \cite{Chikuse:2003aa,Muirhead:1982fk,Gross:1987bf}. It satisfies $C_{\lambda}(X) = 0$ if $\rank(X) < \ell(\lambda)$, and 
\begin{equation}
  \label{eq:intCkappa}
  \int_{\mathcal O} C_{\lambda}(XO YO^\top ) \d\sigma_{\mathcal O(d)}(O) =
  \frac{C_{\lambda}(X) C_{\lambda}(Y)}{C_{\lambda}(I_{d})}, \qquad X,Y \in \R^{d\times d}_{\sym},\quad \ell(\lambda) \le d,
\end{equation}
where $I_{d}$ is the identity matrix in $\R^{d}$. They are normalized such that $
  \Tr(X)^{t} = \sum_{
      |\lambda|=t 
  } C_{\lambda}(X)$, for $X \in \R^{d\times d}_{\sym}.
$ 
The evaluation of zonal polynomials can be extended, 
so that $C_\lambda(XY)$ makes sense for $X,Y\in \R^{d\times d}_{\sym}$,
cf.~\cite{Gross:1987bf}. 
The zonal polynomials do not depend on the dimension
$d$ of the matrix argument, 
but only on the nonzero eigenvalues of the
matrix, which explains the lack of an index $d$ in their notation. 

%

\subsection{Intertwining functions}\label{sec:intertw func}
Let $\lambda$ be a partition with $\ell(\lambda)\leq d/2$. For $P \in \G_{k,d}$, $Q \in \G_{{l},d}$ with $\ell(\lambda)\leq k,{l}\leq d-\ell(\lambda)$, the intertwining functions $p_{\lambda}^{k,{l}}$ considered in Section \ref{sec:L2kernels}, cf.~\eqref{eq:normalizedplk} and \eqref{eq:pkl_rel}, satisfy the symmetry relations 
\begin{equation*}
p_{\lambda}^{k,{l}}(P,Q) = p_{\lambda}^{{l},k}(Q,P),\qquad p_{\lambda}^{k,{l}}(P,Q) = \varepsilon_{k} p_{\lambda}^{d-k,{l}}(I-P,Q)
\end{equation*}
where the constant $\varepsilon_{k} \in \{ \pm 1\}$ only depends on $k$. The latter is a consequence of $P\mapsto I-P$ inducing an intertwining operator.

\subsection{Jacobi polynomials}\label{sec:Jacobi}
In \cite{James:1974aa} the intertwining functions $p_{\lambda}^{k,{l}}$ are related
to the generalized Jacobi polynomials $J_\lambda^{\frac {l} 2,\frac d 2}$ in $k$
variables only for the cases $1 \le k \le {l} \le \tfrac d 2$, $\ell(\lambda)\leq k$ via \eqref{eq:pkl_jacobi}, where 
\begin{equation*}
J_{\lambda}^{\frac\ell 2,\frac d 2}(y_{1}(PQ),\dots,y_{k}(PQ)) = \sum_{\lambda'\le \lambda}
  c_{\lambda,\lambda'}^{\frac d 2} q_{\lambda,\lambda'}(\tfrac{k}{2}) q_{\lambda,\lambda'}(\tfrac {l} 2)
  C_{\lambda'}(PQ),
\end{equation*}
and $y_{1}(PQ) \ge \dots \ge y_{d}(PQ)$ are the eigenvalues of the
  matrix $P Q$ (counted with multiplicities). In order to deal with general $k,{l}$, we define generalized Jacobi polynomials $J_{\lambda}^{\alpha,\beta}:\R^{m} \to \R$ beyond the usual parameter range $\frac12(m-1) < \alpha < \beta - \frac12(m-1)$ in \cite{Dumitriu:2007ax,James:1974aa,Davis:1999zg,Davis:1999dn}, i.e., we define, for any $\beta > \ell(\lambda) - 1$, $\alpha \in \R$, and $y_{1},\dots,y_{m} \in \R$, $m \in \N$, 
\begin{equation}
  \label{eq:jacobi_ext}
  J_{\lambda}^{\alpha,\beta}(y_{1},\dots,y_{m}):=\sum_{\lambda'\le \lambda}
  c_{\lambda,\lambda'}^{\beta} q_{\lambda,\lambda'}(\alpha) q_{\lambda,\lambda'}(\tfrac m 2)
  C_{\lambda'}(\diag(y_{1},\dots,y_{m})).
\end{equation}
The functions $q_{\lambda,\lambda'}$ are continuously and uniquely extended to $\R$. The coefficients 
$c_{\lambda,\lambda'}^{\beta}$ in \eqref{eq:jacobi_ext} depend on $\lambda,\lambda',\beta$
in a rational way, due to the recursion formula in \cite{Dumitriu:2007ax,James:1974aa}. In particular, they are well-defined for $\beta > \ell(\lambda) - 1$. 

Certain properties and symmetry relations of the generalized Jacobi polynomials in the usual parameter range are also valid in the extended parameter range. Indeed, let $\lambda$ be a partition satisfying $\beta > \ell(\lambda) - 1$, $\alpha \in \R$. For $m\in\N$, one observes
\begin{equation*}      J_{\lambda}^{\alpha,\beta}(y_{1},\dots,y_{m})  = (-1)^{\lambda}
      J_{\lambda}^{\beta-\alpha,\beta}(1- y_{1},\dots,1-y_{m}),\qquad y_{1},\dots,y_{m} \in \R.
      \end{equation*}
If $m < \ell(\lambda)$, then $
  J_{\lambda}^{\alpha,\beta}(y_{1},\dots,y_{m})=0$ holds. If $k,{l} \in \N$ with $ {l} > k$, then we obtain
\begin{equation*}
 J_{\lambda}^{\tfrac {l} 2,\beta}(y_{1},\dots,y_{k}) = J_{\lambda}^{\tfrac k
      2,\beta}(y_{1},\dots,y_{k},\overbrace{0,\dots,0}^{{l} -
      k\text{\normalfont -times}}),\qquad y_{1},\dots,y_{k} \in \R.
      \end{equation*}


\subsection{Even and odd kernels}
The sets of even and odd functions on $\G_{d}$ are denoted by 
  \begin{equation*}
    S_{+}  := \{ f \in L^{2}(\G_{d}) \;|\; f=f(I-\cdot)\}, \qquad
    S_{-}  := \{ f \in L^{2}(\G_{d}) \;|\; f=-f(I-\cdot) \},
  \end{equation*}
respectively, and the orthogonal decomposition
$L^{2}(\G_{d}) = S_{+} \oplus S_{-}$ holds. 
We call positive definite zonal kernels $K_{+}$ and $K_{-}$ even and odd if
they satisfy
\[
  K_{+}(P,Q)=K_{+}(I-P,Q), \quad K_{-}(P,Q)=-K_{-}(I-P,Q), \qquad P,Q \in
  \G_{d},
\]
respectively. The corresponding reproducing kernel Hilbert spaces
$\mathcal S(K_{\pm})$ contain only odd and even functions, respectively. If $K$ satisfies $K(P,Q) = K(I-P,I-Q)$, for $P,Q \in \G_{d}$, then it can be uniquely decomposed into an even and odd part, i.e., $K = K_{+} + K_{-}$ such that
  \begin{equation*}
    \mathcal S(K) = \mathcal S(K_{+}) \oplus \mathcal S(K_{-}),\quad \mathcal
    S(K) \cap S_{+} = \mathcal S(K_{+}),\quad \mathcal S(K) \cap S_{-} =
    \mathcal S(K_{-}).
  \end{equation*}
The kernel $K^{s}(P,Q) := 1 + 2\Tr(PQ)^{s} + 2\Tr((I-P)(I-Q))^{s}$, for $P,Q \in \G_{d}$, linearly generates $\Pol_s(\G_d)$, cf.~Remark~\ref{rem:erster}. One observes that $K_{+}^{2t+1}(P,\cdot)$ is a polynomial of degree $2t$ and $K_{-}^{2t+2}(P,\cdot)$ is a polynomial of degree $2t+1$, cf.~the proof of Proposition~\ref{prop:2t to 2t+1}. Therefore, we obtain
  \begin{equation}
    \label{eq:poly_even_odd}
      \Pol_{2t+1}(\G_{d}) \cap S_{+}  \subset \Pol_{2t}(\G_{d}),\qquad 
      \Pol_{2t+2}(\G_{d}) \cap S_{-}  \subset  \Pol_{2t+1}(\G_{d}).
  \end{equation}

Hence, since $H^{|\lambda|}_\lambda(\mathcal{G}_{d})$ has multiplicity one, its reproducing kernel $p_\lambda$ is an even or odd kernel, depending on $|\lambda|$, i.e., it satisfies the symmetry relations
 \begin{equation}
     \label{eq:ppi_sym}
     p_{\lambda}(P,Q) = (-1)^{|\lambda|} p_{\lambda}(P,I-Q), \qquad P,Q \in \G_{d},\qquad \ell(\lambda)\leq d/2.
   \end{equation}

We shall verify that being even or odd transfers into symmetry conditions on the kernel's coefficients when expanded in zonal polynomials.
\begin{lemma}
  \label{lem:r_pi}
  For any partition $\lambda$ with $\ell(\lambda) \le \tfrac d2$, let
  $K:\G_{d}\times \G_{d} \to \R$ be a positive definite zonal kernel obeying $ K(P,Q)=(-1)^{t} K(P,I-Q)$, for $P,Q \in \G_{d}$, 
  with some $t \ge |\lambda|$, and admitting the expansion
\begin{equation}\label{eq:com aa}
      K(P,Q) = \sum_{|\lambda'| \le t, \; \ell(\lambda') \le d}
      r_{\lambda'}(\Tr(P),\Tr(Q)) C_{\lambda'}(PQ), \qquad P,Q \in \G_{d}.
\end{equation}
  If $\mathcal S(K) \subset H_{\lambda}(\G_{d})$, then 
  it holds, for
  $k \in \{\ell(\lambda), \dots,\tfrac d 2\}$ and ${l}\in \{\ell(\lambda),\dots,d-\ell(\lambda)\}$,
  \begin{equation}
    \label{eq:r_pi}
   r_{\lambda}(k,{l}) =
    (-1)^{t-|\lambda|}r_{\lambda}(k,d-{l}).
  \end{equation}
\end{lemma}
The following proof is based on comparing coefficients in the expansion \eqref{eq:com aa} after exploiting the binomial expansion of zonal polynomials.
\begin{proof}
  We know by Corollary~\ref{cor:basischange}
  that the only nonzero Fourier coefficient of $K$ is $\widehat{K}_{\lambda}$. 
We now fix $k\in \mathcal{I}:=\{\ell(\lambda),\dots,\lfloor \tfrac d 2 \rfloor\}$. Comparing coefficients of 
  $C_{\lambda'}(PQ)$ on $\G_{k,d}\times \G_{{l},d}$, ${l} \in \mathcal{I}$ via \eqref{eq:pkl_jacobi} yields 
$
    r_{\lambda'}(k,{l}) = 0$, for ${l} \in \mathcal{I}$, $|\lambda'|>|\lambda|$, $\ell(\lambda')\leq k$. 
  For $P\in\G_{k,d}$, the orthogonal invariance and the binomial expansion of the zonal polynomials, cf.~\cite{Dumitriu:2007ax,Muirhead:1982fk,Chikuse:2003aa}, yield
  \[
    \begin{aligned}
      K(P,I-Q) & = \sum_{|\lambda'| \le t,\;\; \ell(\lambda') \le k}
      r_{\lambda'}(\Tr(P),d-\Tr(Q)) C_{\lambda'}(P(I-Q))\\
      & = \sum_{|\lambda'| \le t,\;\; \ell(\lambda') \le k}
      r_{\lambda'}(\Tr(P),d-\Tr(Q)) \sum_{\lambda'' \le \lambda'}\binom{\lambda'}{\lambda''}(-1)^{|\lambda''|}\frac{C_{\lambda'}(I_{k})}{C_{\lambda''}(I_{k})}C_{\lambda''}(PQ)\\
      & = \sum_{|\lambda''| \le t,\;\; \ell(\lambda'') \le k} \tilde
      r_{\lambda''}(\Tr(P),\Tr(Q)) C_{\lambda''}(PQ),\qquad P,Q \in \G_{d},\quad
    \end{aligned}
  \]
  where $\binom{\lambda'}{\lambda''}$ are generalized binomial coefficients, cf.~\cite{Dumitriu:2007ax}, and, for $|\lambda''|\le t$ with $\ell(\lambda'')\le k$, 
  \begin{equation}
    \label{eq:r_pi_binom}
    \tilde r_{\lambda''}(k,{l}) := \sum_{\substack{|\lambda'| \le t,\\\lambda'' \le \lambda',\; \ell(\lambda')\leq k}}
    \binom{\lambda'}{\lambda''}(-1)^{|\lambda''|}\frac{C_{\lambda'}(I_{k})}{C_{\lambda''}(I_{k})}
    r_{\lambda'}(k,d-{l}), \qquad {l} \in \{1,\dots,d-1\}.
  \end{equation}
According to the assumption $ K(P,Q)=(-1)^{t} K(P,I-Q)$, we again compare coefficients and obtain 
$
    \tilde r_{\lambda''}(k,{l}) = 0$, for ${l} \in \mathcal{I}$, $ |\lambda''| > |\lambda|$, $ \ell(\lambda'')\leq k$.
 Note that $\binom{\lambda'}{\lambda'}=1$ holds, cf.~\cite{Dumitriu:2007ax}. One can now show by induction using \eqref{eq:r_pi_binom} and starting with
  $|\lambda''|=t$ that this implies 
$
    r_{\lambda'}(k,d-{l}) = 0$, ${l} \in \mathcal{I}$, $|\lambda'| > |\lambda|$, $ \ell(\lambda')\leq k$. 
 By using the latter, for $P\in\G_{k,d}$ and $Q \in \G_{\mathcal{I},d}$, we
  arrive  at the expansions
    \begin{align*}
      K(P,Q) = & \hspace{-2ex}\sum_{|\lambda'| \le |\lambda|,\;\; \ell(\lambda') \le k}
     \hspace{-2ex} r_{\lambda'}(\Tr(P),\Tr(Q)) C_{\lambda'}(PQ),\\
      K(P,I-Q) 
       = &\hspace{-4ex}\sum_{|\lambda''|=|\lambda|,\;\;\ell(\lambda'')\le k}      \hspace{-5ex} (-1)^{|\lambda''|}
r_{\lambda''}(\Tr(P),d-\Tr(Q)) C_{\lambda''}(PQ) 
       + \hspace{-4ex}\sum_{|\lambda''|<|\lambda|,\;\;\ell(\lambda'')\le k} \hspace{-5ex}\tilde
      r_{\lambda''}(\Tr(P),\Tr(Q)) C_{\lambda''}(PQ).
    \end{align*}
According to $K(P,Q)=(-1)^{t} K(P,I-Q)$ the assertion 
  \eqref{eq:r_pi} is derived by comparing coefficients again.
  \end{proof}

\section{Proof of Theorem~\ref{the:ki}}\label{app:proof of long}
Let us complete the proof of Theorem~\ref{the:ki} by verifying that \eqref{eq:Ski_inc} holds with equality. 
\begin{proof}[Proof of Theorem~\ref{the:ki}]
%
We shall verify equality in \eqref{eq:Ski_inc} by induction over $s$. For $s=0$, we observe
  $\mathcal S(K_{\lambda}^{0}) = \mathcal S(p_{\lambda}) = H_{\lambda}^{|\lambda|}(\G_{d})$. Now, fix $0 < s \le |\{\ell(\lambda),\ldots,d-\ell(\lambda)\}|-1$ and assume that the equality holds for all $s'$ with $0 \le s' < s$. We consider an arbitrary
  positive definite zonal kernel
  $K:\R^{d\times d}_{\sym} \times \R^{d\times d}_{\sym} \to \R$ of degree at
  most $|\lambda|+s$ such that
  \begin{equation}
    \label{eq:KHpis}
   \mathcal S(K|_{\G_{d}\times
      \G_{d}}) \cong \mathcal H_{2\lambda}^{d},\qquad  \mathcal S(K|_{\G_{d}\times \G_{d}}) \subset H_{\lambda}^{|\lambda|+s}(\G_{d}),\qquad  \mathcal S(K|_{\G_{d}\times \G_{d}})\perp
H_{\lambda}^{|\lambda|+s-1}(\G_{d}) .
  \end{equation}
According to invariant theory, cf.~\cite[Theorem~7.1]{Procesi:1976ul}, the restriction $K|_{\G_{d}\times\G_{d}}$ is a linear combination of terms of the form
  \[
    \Tr(X)^{k'}\Tr(Y)^{{l}'} \prod_{i=1}^{\ell(\lambda')}\Tr((XY)^{\lambda'_{i}}),\qquad
    k'+|\lambda'|,\; {l}'+|\lambda'| \le |\lambda|+s.
  \]
Since any polynomial $\prod_{i=1}^{\ell(\lambda')} \Tr(X^{\lambda'_{i}})$ can be
  expanded into zonal polynomials $C_{\lambda''}(X)$ with $|\lambda''| = |\lambda'|$, there are polynomials $r_{\lambda'}:\R\times\R\to\R$ of degree at most
  $|\lambda|+s-|\lambda'|$, such that 
  \begin{equation*}
      K(P,Q) = \sum_{\substack{|\lambda'| \le |\lambda|+s}}
      r_{\lambda'}(\Tr(P),\Tr(Q)) C_{\lambda'}(PQ), \qquad P,Q \in \G_{d}.
  \end{equation*}
Moreover, due to \eqref{eq:poly_even_odd}
  and the second and third 
  property in \eqref{eq:KHpis}, the kernel $K|_{\G_{d}\times \G_{d}}$ is also
  even or odd if $|\lambda|+s$ is even or odd, respectively, i.e.,
  \begin{align}
    K(P,Q) &= (-1)^{|\lambda|+s}K(P,I-Q), \qquad P,Q \in \G_{d}.    \label{eq:K_sym}
\intertext{According to Corollary~\ref{cor:basischange} and $ \mathcal S(K|_{\G_{d}\times
      \G_{d}}) \cong \mathcal H_{2\lambda}^{d}$, there is a function $q:\R \to \R$, such that }
    K(P,Q) &= q(\Tr(P))q(\Tr(Q)) p_{\lambda}(P,Q), \qquad P,Q \in
    \G_{d}.\nonumber
  \end{align}
We now compare coefficients of $C_{\lambda}(PQ)$ on
  $\G_{k,d}\times \G_{{l},d}$,
  $k,{l} \in \mathcal{I}:=\{\ell(\lambda),\dots,\lfloor \tfrac d 2 \rfloor\}$ in
  both expansions. We obtain $ r_{\lambda}(k,{l}) = q(k)q({l})$, for
  $k,{l} \in \mathcal{I}$. Since 
  $ \mathcal S(K|_{\G_{d}\times \G_{d}}) \cong \mathcal H_{2\lambda}^{d}$, the kernel $K|_{\G_{d}\times \G_{d}}$ does not completely vanish on $\G_{d}$, so that there is $k' \in \mathcal{I}$ such
  that $q(k')\ne 0$ implying  $q({l}) = r_{\lambda}(k',{l})/q(k')$, for
  ${l} \in \mathcal{I}$. Furthermore, the symmetry \eqref{eq:K_sym} and
  \eqref{eq:ppi_sym} implies $q(d-l) = (-1)^{s}q(l)$ for ${l} \in\mathcal{I}$, so that Lemma~\ref{lem:r_pi} yields
    \[
    q(d-{l}) = (-1)^{s}q({l}) = (-1)^{s} r_{\lambda}(k',{l})/q(k') = 
    r_{\lambda}(k',d-{l})/q(k'), \qquad {l} \in\mathcal{I}.
  \]
  Hence, $q(l) = r_\lambda(k',l)/q(k')$,
  $l \in \{\ell(\lambda), \dots, d-\ell(\lambda) \}$ and 
  $q|_{\{\ell(\lambda),\dots,d-\ell(\lambda)\}}$ is the restriction of a polynomial of degree at
  most $s$. The choice of $K$ and the induction hypothesis yield 
  $K|_{\G_{d}\times\G_{d}} * K_{\lambda}^{i} \equiv 0$, for $i=0,\dots,s-1$,
  which is equivalent to the orthogonality conditions 
  \[
    \sum_{m=\ell(\lambda)}^{d-\ell(\lambda)} q_{\lambda}^{i}(m) q(m) |v_{\lambda}^{m}|^{2}=0, \qquad
    i=0,\dots,s-1.
  \]
  Thus, $q_{\lambda}^{s}|_{\{\ell(\lambda),\dots,d-\ell(\lambda)\}}$ and $q|_{\{\ell(\lambda),\dots,d-\ell(\lambda)\}}$ are linearly dependent, so that $K|_{\G_{d}\times \G_{d}}=c K_{\lambda}^{s}$, for some constant $c\neq 0$, 
implying 
  $\mathcal S(K|_{\G_{d}\times \G_{d}}) = \mathcal S(K_{\lambda}^{s})$. Since $K$
  was an arbitrary kernel with the property \eqref{eq:KHpis}, we derive
$ 
    H_{\lambda}^{|\lambda|+s}(\G_{d})  / H_{\lambda}^{|\lambda|+s-1}(\G_{d}) \cong
    \mathcal H_{2\lambda}^{d},
  $ 
  so that \eqref{eq:Ski_inc} yields the assertion for $s$.
\end{proof}

\section{The reproducing kernels for $|\lambda|=0,1,2$}\label{app:ex}
Let $K_{2\lambda}:\R^{d\times d}_{\sym} \times \R^{d\times d}_{\sym} \to \R$ denote the reproducing kernel (with respect to the differentiation inner product) of the irreducible representation $\mathcal H_{2\lambda}^{d}$ in $\Hom_{|\lambda|}(\R^{d\times d}_{\sym})$. For $|\lambda|=0,1,2$ and $X, Y \in \R^{d\times d}_{\sym}$, we obtain
{\small
\begin{align*}
  K_{(0)}(X,Y) & = 1,\qquad\qquad  K_{(2)}(X,Y) = \Tr(X Y) - \frac{\Tr(X)\Tr(Y)}{d},\\
    K_{(4)}(X,Y) &= \frac{\Tr(XY)^2+2\Tr((XY)^2)}{6} +
              \frac{\big(\Tr(X)^2+2\Tr(X^2)\big)\big(\Tr(Y)^2+2\Tr(Y^2)\big)}{6(d+2)(d+4)}\\
           & - \frac{\Tr(XY)\Tr(X)\Tr(Y)+2\Tr(XY^2)\Tr(X)+2\Tr(X^2Y)\Tr(Y)+4\Tr(X^2Y^2)}{3(d+4)},\\
  K_{(2,2)}(X,Y)&=\frac{\big(\Tr(XY)^2-\Tr((XY)^2)\big)}{3}  +
                \frac{\big(\Tr(X)^2-\Tr(X^2)\big)\big(\Tr(Y)^2-\Tr(Y^2)\big)}{3(d-1)(d-2)}\\
           &\quad - \frac{\Tr(XY)\Tr(X)\Tr(Y)-\Tr(XY^2)\Tr(X)-\Tr(X^2Y)\Tr(Y)+\Tr(X^2Y^2)}{3(d-2)/2}.
\end{align*}
}
The kernel $K_{(2,2)}$ is only
defined for $d \ge 3$ and $K_{(2,2)}\equiv 0$ for $d=3$. This is
in accordance with the decomposition of
$\Hom_{2}(\R^{d\times d}_{\sym})$, where the irreducible representation
$\mathcal H_{(2,2)}^{d}$ only occurs for $d \ge 4$, cf.~\eqref{eq:branch rules explit}. The restriction $K_{2\lambda}|_{\G_{d}\times \G_{d}}$ coincides with $p_\lambda$ up to a multiplicative constant. Indeed, we obtain
  \begin{align*}
p_{(0)}& = 1|_{\mathcal{G}_{d}\times \mathcal{G}_{d}}, &
p_{(2)}& = 2 K_{(4)}|_{\mathcal{G}_{d}\times \mathcal{G}_{d}},\\
p_{(1)}&= K_{(2)}|_{\mathcal{G}_{d}\times \mathcal{G}_{d}},&
p_{(1,1)}&= 2 K_{(2,2)}|_{\mathcal{G}_{d}\times \mathcal{G}_{d}}.
\end{align*}
Further restricting $p_\lambda$ to $\mathcal{G}_{\mathcal{I},d}\times \mathcal{G}_{\mathcal{I},d}$ yields the reproducing kernel of $H^{|\lambda|}_\lambda(\G_{\mathcal{I},d})$. 
The intertwining polynomials $p_\lambda^{k,{l}}=(v_{\lambda}^{k}v_{\lambda}^{{l}})^{-1}p_\lambda|_{\G_{k,d}\times \G_{{l},d}}$ are derived from 
\begin{align*}
v_{(0)}^{k}& = 1, & 
v_{(2)}^{k}& =  \sqrt{\frac{8k(k+2)(d-k)(d-k+2)}{(d-1)d(d+1)(d+2)(d+4)(d+6)}},\\
v_{(1)}^{k} & = \sqrt{\frac{2k(d-k)}{(d-1)d(d+2)}} , & 
v_{(1,1)}^{k} & = \sqrt{\frac{8(k-1)k(d-k-1)(d-k)}{(d-3)(d-2)(d-1)d(d+1)(d+2)}}.
\end{align*}
Reproducing kernels of the subspaces of $H^{|\lambda|+s}_\lambda(\G_{\mathcal{I},d})$ equivalent to $\mathcal{H}^d_{2\lambda}$ with increasing polynomial degree can be derived by restrictions of the kernels $K^i_\lambda$ in the proof of Theorem~\ref{the:ki}, cf.~\eqref{eq:Ki}.

\end{appendix}


\end{document}